\documentclass[11pt, a4paper]{article}
\usepackage{latexsym}
\usepackage{indentfirst}
\usepackage{graphicx}
\usepackage{amsmath}
\usepackage{amssymb}
\begin{document}
\title{Study of the behaviour of proliferating cells in leukemia  modelled by a system of delay differential equations}
\author{Anca Veronica Ion\\"Gh. Mihoc-C. Iacob" Institute of Mathematical Statistics\\
and Applied Mathematics of the Romanian Academy, \\13, Calea 13 Septembrie, 050711, Bucharest, Romania;\\ email: anca-veronica.ion@ima.ro.} \date{}
\maketitle

\begin{abstract}
For the model of periodic chronic myelogenous leukemia considered by
Pujo-Menjouet, Mackey et al., model consisting of two delay
differential equations, the equation for the density of so-called
"resting cells" was studied from numerical and qualitative point of
view in several works. In this paper we focus on the equation for
the density of proliferating cells and study it from a qualitative
point of view.

\vspace{0.2cm}

\end{abstract}

\section{Introduction}

We consider the model of periodic chronic myelogenous leukemia consisting in the system of two delay differential equations \cite{PM-M}, \cite{PM-B-M}:

\[\quad \frac{dP(t)}{dt}=-\gamma P(t) +\beta(N(t))N(t)-
\frac{k}{2}\beta(N(t-r))N(t-r),
\]
\[\frac{dN(t)}{dt}=-\left[\beta(N)+\delta
\right]N(t)+k\beta(N(t-r))N(t-r).
\]
Here $P(t)$ represents the density of proliferating cells while $N(t)$ is the density of so-called ``resting" cells.

The function $\beta$ is given by $\displaystyle \beta(N)=\beta_0\theta^n/(\theta^n+N^n),\,n>0,$ and represents the rate at which the resting cells re-enter the proliferating phase.

\vspace{0.3cm}

For the detailed description of the model see \cite{PM-M}. We just remind that $\gamma>0$ is the rate of death of proliferating cell by apoptosis,  $k=2e^{-\gamma r},$ $\delta>0$ is the rate by which the resting cells transform into normal cells, $\beta_0>0$ is the maximal rate of cell movement from the resting
phase into proliferation; $\theta$ is the density of resting cell population such that $\beta(\theta)=\beta_0/2.$ The delay $r$ refers to a processus of division of the proliferation cells that takes place at time $r$ after they entered into the proliferating phase \cite{PM-M}.

\vspace{0.3cm}

With $x(t)=P(t)/\theta,\,\,y(t)=N(t)/\theta$ the system becomes:

\begin{equation}\label{eq-x}\;\;\;\;\;\;\;\dot{x}(t)=-\gamma x(t) +\frac{\beta_0y(t)}{1+y(t)^n}-
\frac{k}{2}\frac{\beta_0y(t-r)}{1+y(t-r)^n},
\end{equation}
\begin{equation}\label{eq-y}\dot{y}(t)=-\left[\frac{\beta_0}{1+y(t)^n}+\delta\right]y(t)+k\frac{\beta_0y(t-r)}{1+y(t-r)^n}.
\end{equation}

\vspace{0.2cm}

We see that the equation for $y$ does not depend on $x$. It was
studied, from a numerical point of view in \cite{PM-M},
\cite{PM-B-M}, and from a qualitative point of view in
\cite{AI-stab}, \cite{AI-Hopf}, \cite{AVI-RMG} (in all these works,
the density of resting cells,denoted here by $y$, was denoted by
$x$). More precisely, in \cite{AI-stab}, the authors studied in
depth the stability of the two equilibria of equation \eqref{eq-y},
in \cite{AI-Hopf} the stability of periodic solution emerged by Hopf
bifurcation was considered, while in \cite{AVI-RMG} the existence of
Bautin type bifurcation points was established.

Here we consider the first equation, that for the proliferating
cells. Since their increase is the cause of the illness, the
physician should be interested in the behaviour of $x(t)$ more than
in that of $y(t)$. The equation for $x$ is very simple, it is a
linear forced equation, with the forcing term depending on $y$.
Hence the behaviour of $x$ is determined by that of $y$ (in the
absence of the terms containing $y$, $x(t)$ would tend to zero when
$t\rightarrow \infty$).

\vspace{0.3cm}

In the study of eq. \eqref{eq-y}, (\cite{AI-stab}, \cite{AI-Hopf}, \cite{AVI-RMG}) in the zone of parameters that we considered (see Section 4) we encountered the following typical behaviours of $y$:

\begin{itemize}
  \item $y(\cdot)$ is an equilibrium solution of \eqref{eq-y}, i.e. $y(t)=ct$;
  \item $y(t)$ tends to an equilibrium solution of \eqref{eq-y}, when $t$ tends to $\infty$;
  \item $y(\cdot)$ is a periodic solution of \eqref{eq-y};
  \item $y(t)$ tends to a periodic solution of \eqref{eq-y}, when $t$ tends to $\infty$.
\end{itemize}

In the present paper, by using the results obtained on $y(\cdot)$, we study the behavior of $x(\cdot)$.

In the Appendix we prove some propositions that show how the behaviour of $y$ influences that of $x$. These are used in all the main sections of the paper.

In Section 2 we discuss the stability of the equilibria of eq. \eqref{eq-x}, indicating the zone of the parameters space where these are stable. Section 3 refers to the occurrence of Hopf bifurcation for $y$, that induces a similar behaviour for $x$. Section 4 refers to the Bautin bifurcation of $y$, that was studied in \cite{AVI-RMG} and to its influence on $x$.

A section of Conclusions summarizes all the results in the paper.

\section{Equilibria of the system, their stability}

The system \eqref{eq-x} - \eqref{eq-y} has the equilibria:
$x_1=0,\;\;y_1=0,$ and $$x_2=\frac{2-k}{2\gamma}\frac{\beta_0
y_2}{1+y_2^n},\;\;y_2=\left(\frac{\beta_0}{\delta}(k-1)-1\right)^{1/n}.$$

In this section we study the stability of these equilibria, by using the stability results contained in \cite{AI-stab}.
Let us remind that $y(\cdot),$ the solution of \eqref{eq-y}, is a bounded function, defined on $[-r,\,\infty)$, as shown in \cite{AI-stab}.

\vspace{0.3cm}

\subsection{Stability of ($x_1$, $y_1$)}

We intend to establish the stability of $(x_1,\,y_1)=(0,\,0)$ depending on the parameters of the problem. For this we shall use the results on the stability of $y$ in \cite{AI-stab} and the results in the Appendix.

\vspace{0.3cm}

 \textbf{Propopsition 2.1.}   \textit{If \begin{equation}\label{cond-stab-0} \frac{\beta_0}{\delta}(k-1)<1,\end{equation} then the equilibrium solution $(x_1,\,y_1)$ of system \eqref{eq-x}, \eqref{eq-y} is   asymptotically stable.}

\vspace{0.2cm}

\textbf{Proof.}  By determining the set of parameters for which all eigenvalues of the linear equation attached to the delay differential equation  \eqref{eq-y} have negative real parts, in \cite{AI-stab} it is shown that the solution $y_1=0$ of eq. \eqref{eq-y} is  asymptotically stable when $\frac{\beta_0}{\delta}(k-1)<1$.

Hence, in the same set of parameters, if the initial function
$\phi\in \mathcal{C}([-r,0],\,\mathbb{R})$ for equation \eqref{eq-y}
is close enough to the function $y_1(s)=0,\,s\in[-r,0]$ such that
the corresponding solution, $y(t,\,\phi)$ tends to $0$ when
$t\rightarrow \infty,$ Proposition A.1 of the Appendix implies that
$\displaystyle \lim_{t\rightarrow \infty}x(t)= 0.$

 It follows that the equilibrium solution $(x_1,\,y_1)$ of system \eqref{eq-x}, \eqref{eq-y} is asymptotically stable.$\Box$

\vspace{0.3cm}

If we think at eq. \eqref{eq-x} as a nonautonomous equation (and not a part of an autonomous system), in the conditions of the above Proposition, and if $\phi\in \mathcal{C}([-r,0],\,\mathbb{R})$ is such that $\lim_{t\rightarrow \infty}y(t,\phi)= 0,$ then the solution $x_1=0$ of  \eqref{eq-x} is globally asymptotically stable.

\vspace{0.3cm}

\textbf{Proposition 2.2.} \textit{For $\frac{\beta_0}{\delta}(k-1)=1,$ the solution $(x_1,\,y_1)$ of system \eqref{eq-x}, \eqref{eq-y} is  stable.}

\vspace{0.3cm}

\textbf{Proof.} When $\frac{\beta_0}{\delta}(k-1)=1,$ there is an
eigenvalue, $\lambda,$ of the linearized eq. in $y$, that has
$\mathrm{Re} \lambda =0$ \cite{AI-stab}. In this situation, also in
\cite{AI-stab}, we proved that there is a Lyapunov function for eq.
\eqref{eq-y}, hence the solution $y_1=0$ is stable.

\vspace{0.3cm}

Now we prove that the solution $(x_1,\,y_1)$ of  the system is stable. In eq. \eqref{eq-x}, we denote the terms containing $y$ by $F(y(t),\,y(t-r)):=\widetilde{F}(y)(t)$.

For this we consider an $\varepsilon >0$ and chose a $\delta_\varepsilon>0$  with the property that $|\phi|_0<\delta_\varepsilon$ implies $|y(t,\, \phi)| < \min\left(\frac{\varepsilon\gamma}{2\beta_0},\,\varepsilon\right)$, for any $t\geq 0$. This is possible due to the results in \cite{AI-stab}.

Then we have:
\[-\frac{\varepsilon\gamma}{2} < \frac{\beta_0y(t)}{1+y(t)^n}-
\frac{k}{2}\frac{\beta_0y(t-r)}{1+y(t-r)^n}<\frac{\varepsilon\gamma}{2},\,\,\mathrm{for}\,\,t\geq T_\varepsilon.
\]
Since the solution of equation \eqref{eq-x} is
\[
x(t,x_0)=x_0e^{-\gamma t}+e^{-\gamma t}\int_0^t e^{\gamma s}\widetilde{F}(y)(s)ds,
\]
relation  \eqref{sol-int} implies, for $t\geq T_\varepsilon,$
\[|x(t,x_0)|\leq |x_0|e^{-\gamma t}+e^{-\gamma t}\int_0^te^{\gamma s}\frac{\varepsilon\gamma}{2}ds=
\]
\[=|x_0|e^{-\gamma t}+(1-e^{-\gamma t})\frac{\varepsilon}{2} < |x_0|+\frac{\varepsilon}{2}.
\]
Now, for $|x_0|\leq \frac{\varepsilon}{2},$ and $t\geq 0,$ we have
\[|x(t,x_0)|\leq \varepsilon,
\]
and since $|y(t,\, \phi)| <\varepsilon$ the solution $(0,\,0)$ of system \eqref{eq-x}, \eqref{eq-y} is stable.$\Box$

\vspace{0.3cm}

Again, if we think at eq. \eqref{eq-x} as a nonautonomous equation, in the conditions of Proposition 2.2, and if the initial condition of \eqref{eq-y}, $\phi\in \mathcal{C}([-r,0],\,\mathbb{R}),$ is close enough to $y_1=0$, then the solution $x_1=0$ of  \eqref{eq-x} is  stable.

\vspace{0.3cm}

\textbf{Proposition 2.3.}  \textit{When $\frac{\beta_0}{\delta}(k-1)>1$,  the solution $(x_1,\,y_1)$ is unstable.}

\vspace{0.2cm}

\textbf{Proof.} The condition $\frac{\beta_0}{\delta}(k-1)>1$ is the condition of existence of the equilibrium $y_2$ of \eqref{eq-y}. In  \cite{AI-stab} we proved that, in this situation, the solution $y_1$ loses stability (becomes unstable), hence the solution $(x_1,\,y_1)$  is also unstable.$\Box$

\subsection{Stability of ($x_2$, $y_2$)}

In the study of stability of $y_2$, \cite{AI-stab}, an important role is played by the constant  $B_1=\frac{d}{dy}\left(\frac{\beta_0y}{1+y^n}\right)\left|_{y=y_2}\right.$ that occurs in the linearized of eq. \eqref{eq-y}. In fact, the linearized around $x_2$ of equation \eqref{eq-y} is
  \begin{equation}\label{lineq}
\dot{z}(t)=-[B_1+\delta]z(t)+kB_1z(t-r),\end{equation}
and the value of $B_1$, in terms of the parameters of the problem is
\begin{equation}\label{B1}
B_1=\frac{\delta}{k-1}\left[\frac{n\delta}{\beta_0(k-1)}-n+1
\right].\end{equation}

\vspace{0.3cm}

\textbf{Proposition 2.4.}\textit{ If the following conditions are fulfilled: }
\begin{itemize}
  \item $B_1<0,$ \;\;($\Leftrightarrow \displaystyle \frac{\beta_0}{\delta}(k-1)>\frac{n}{n-1}$);
  \item $\delta+B_1<0,$ ($\Leftrightarrow n>k\; and \;\displaystyle \frac{\beta_0}{\delta}(k-1)>\frac{n}{n-k}$);
  \item $|\delta+B_1|<|kB_1|$ \textit{and} \;\;$\displaystyle \frac{\arccos{((\delta+B_1)/kB_1)}}{\omega_0}<r<\frac{1}{|\delta+B_1|},$ \textit{where} $\omega_0$ \textit{is the solution in} $(0,\pi/r)$ \textit{of the equation}
$\omega\cot (\omega r)=-(\delta+B_1);$
\end{itemize}
\textit{then, the solution $(x_2,\,y_2)$ of system \eqref{eq-x}, \eqref{eq-y} is asymptotically stable.}

\vspace{0.2cm}
\textbf{Proof.} As shown in \cite{AI-stab}, the conditions in the enounce are  sufficient in order that all the eigenvalues of \eqref{lineq} have negative real part. That is, in these conditions, the solution $y_2$ is asymptotically stable.

For the initial condition $\phi\in \mathcal{C}([-r,0],\mathbb{R})$ in a small neighborhood of $y_2$, such that $\lim_{t\rightarrow \infty} y(t,\phi)=y_2$, Proposition A.1 implies that $\lim_{t\rightarrow \infty} x(t,x_0)=x_2$, for any initial value $x_0$. It follows that $(x_2,\,y_2)$ is asymptotically stable.$\Box$

\vspace{0.3cm}

\textbf{Proposition 2.5.} \textit{If the following conditions are fulfilled:}
\begin{itemize}
  \item $B_1<0$ \textit{and }$\delta+B_1>0,$  ($\Leftrightarrow \displaystyle \frac{n}{n-1}<\frac{\beta_0}{\delta}(k-1)<\frac{n}{n-k}$);
  \item $\delta+B_1>|kB_1|\,\,\mathrm{or}\,\left\{\delta+B_1\leq |kB_1|\,\,
and\,\,\displaystyle r<\frac{\arccos{((\delta+B_1)/kB_1)}}{\omega_0}\right\}$ \textit{with} $\omega_0$ \textit{ defined as above};
\end{itemize}
\textit{then the solution $(x_2,\,y_2)$ of system \eqref{eq-x}, \eqref{eq-y} is asymptotically stable.}

\vspace{0.2cm}

\textbf{Proof.} Again, the results in \cite{AI-stab} show that $y_2$ is asymptotically stable in the conditions of the enounce.
This and Proposition A.1 imply that, for the initial function $\phi$ of eq. \eqref{eq-y} in a convenient neighborhood of $y_2$,  we have $x(t)\rightarrow x_2$ when $t\rightarrow \infty$, hence the equilibrium $x_2$ of equation \eqref{eq-x} is (globally) asymptotically stable, and the equilibrium $(x_2,\,y_2)$ of system \eqref{eq-x}, \eqref{eq-y} is asymptotically stable.$\Box$

\vspace{0.3cm}

If we consider eq. \eqref{eq-x} as a non-autonomous differential equation, if either the conditions of Proposition 2.4 or those of Proposition 2.5 are satisfied and if $\phi$ is such that $\lim_{t\rightarrow \infty} y(t,\phi)=y_2$, then $x_2$ is  a globally asymptotically stable equilibrium of \eqref{eq-x}.

\vspace{0.3cm}

\textbf{Proposition 2.6.}  \textit{If } $B_1>0$, \textit{then the solution $(x_2,\,y_2)$ of system \eqref{eq-x}, \eqref{eq-y} is asymptotically stable.}

\vspace{0.2cm}

\textbf{Proof.} As proved in \cite{AI-stab}, for $B_1>0$, $y_2$ is  asymptotically stable. This implies (via Proposition A.1) the asymptotic stability of $(x_2,\,y_2)$ as solution of eq. \eqref{eq-x}, \eqref{eq-x}.$\Box$

\vspace{0.3cm}

As above, we remark, that, if $\phi$ is such that $\lim_{t\rightarrow \infty} y(t,\phi)=y_2$, then $x_2$ is  a globally asymptotically stable equilibrium of \eqref{eq-x}.

\section{Points of Hopf bifurcation for $y$.\\Influence on $x$}

The analysis in \cite{PM-M},  \cite{PM-B-M}, \cite{AI-Hopf} shows that the equilibrium $y_2$ loses stability by a Hopf bifurcation, at points where $B_1<0$ and
\begin{equation}\label{cond-Hopf}
r=\frac{\arccos{((\delta+B_1)/kB_1)}}{\sqrt{(kB_1)^2-(\delta+B_1)^2}}.
\end{equation}

Propositions A.1, A.2 and A.3 show that, when the delayed differential equation \eqref{eq-y} presents a Hopf bifurcation, the solutions of equation in $x$ behave like having a Hopf bifurcation (because of the term in $y$).

Indeed, for the zone of stability of $y_2$, we already saw that $x(t,\,x_0)\rightarrow x_2,$ when $t\rightarrow \infty$, for any $x_0\geq 0$.

After the Hopf bifurcation, an attractive periodic orbit occurs for the equation in $y$. This implies, as Propositions 2.2 and 2.3 show, that there is a certain initial value $\widetilde{x}_0=\widetilde{u}_0+x_2$ such that $x(t,\widetilde{x}_0)$ is a periodic orbit in the phase portrait of eq. \eqref{eq-x}, orbit that is asymptotically stable.

Hence, in the points from the parameters space where $y(\cdot)$
presents Hopf bifurcation,  the function describing the density of
proliferating cells, $x(\cdot)$, behaves like having a Hopf
bifurcation also. Of course it is not a genuine Hopf bifurcation,
since the eigenvalue of the linearized equation in $x$  at $x_2$ is
$-\gamma$, it is only a behaviour induced by the ``forcing term''
$\widetilde{F}(y)$.

\begin{figure}\centering
\includegraphics[width=0.44\linewidth]{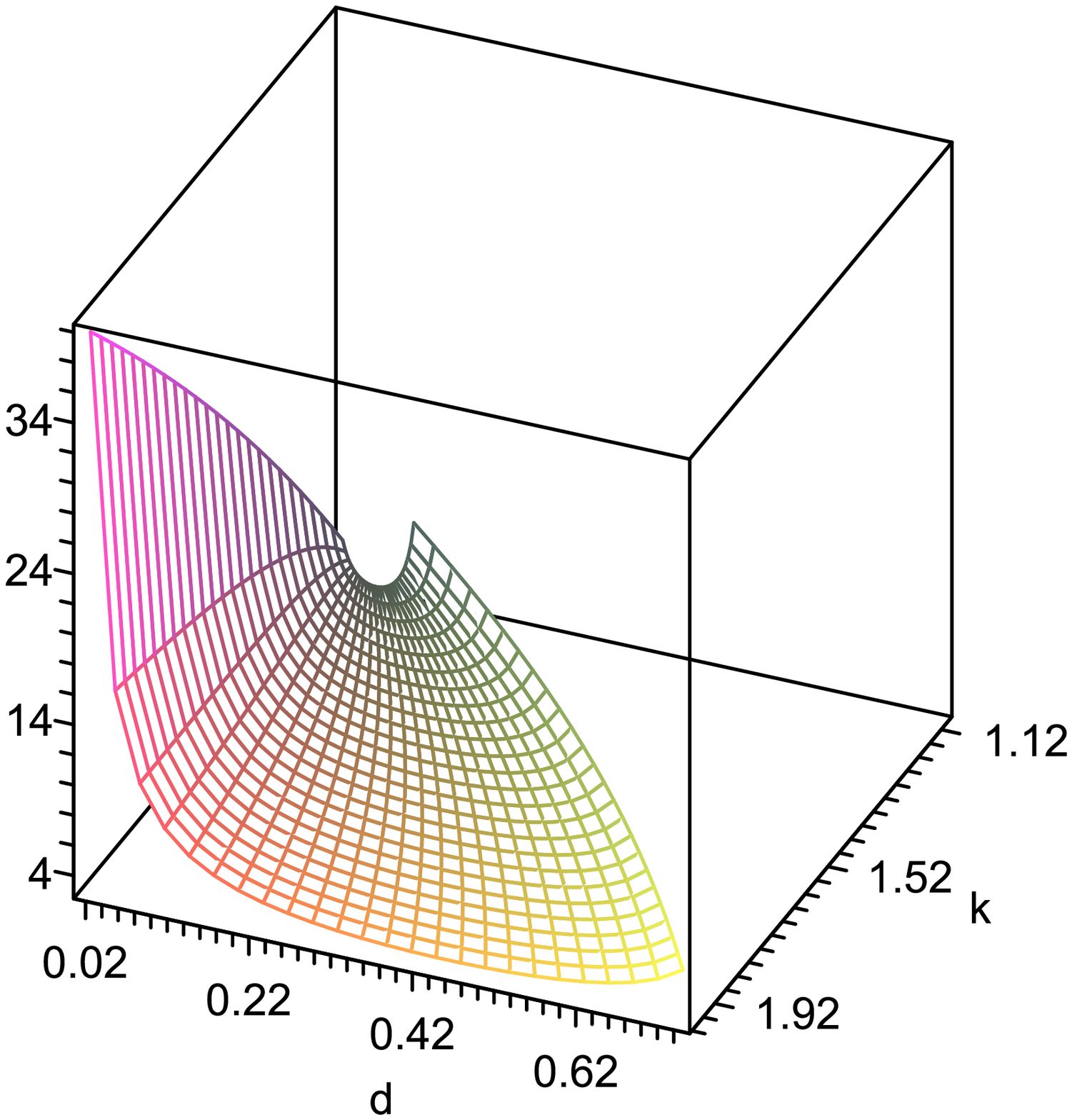}\includegraphics[width=0.49\linewidth]{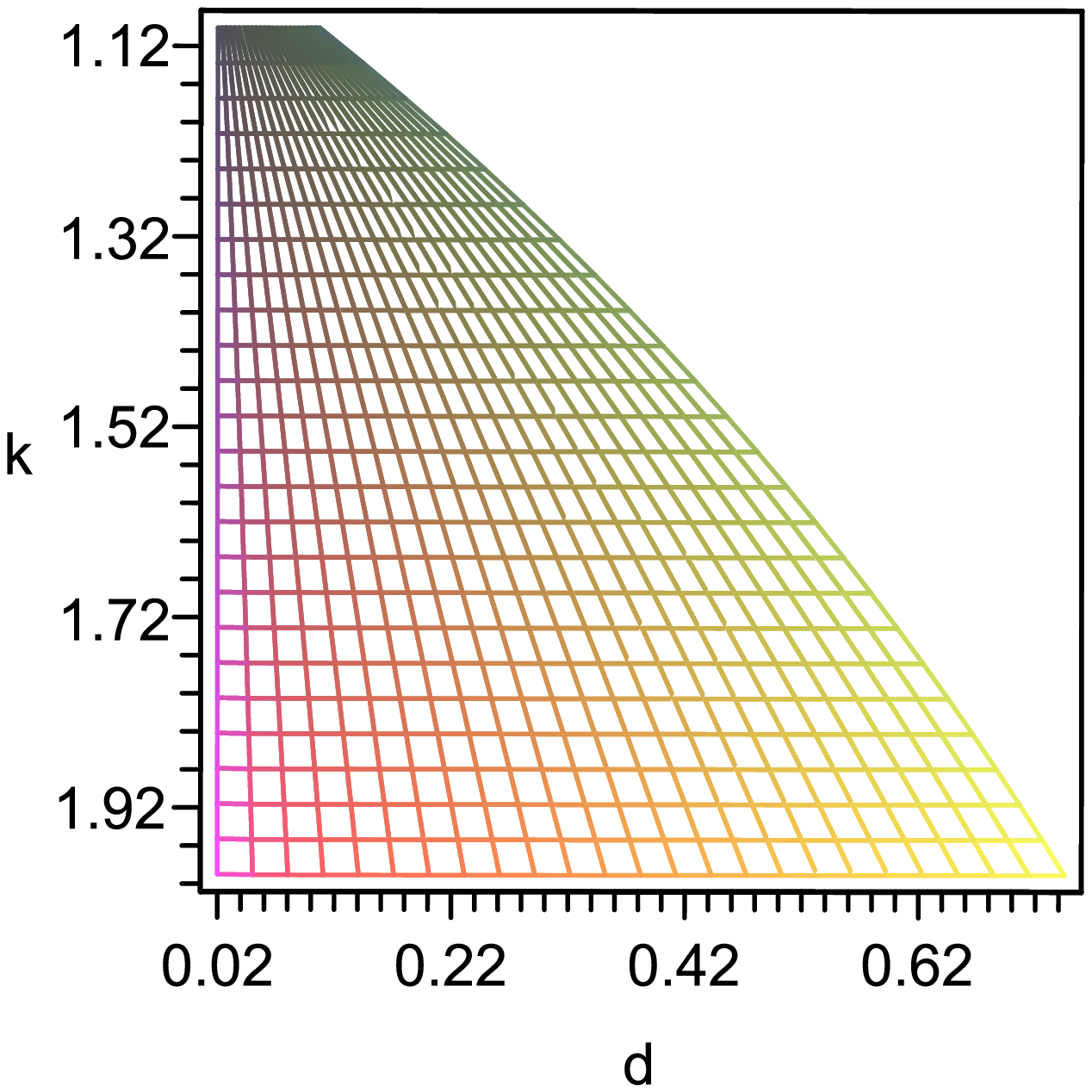}\caption{\small{\textbf{Left}- Surface of Hopf bifurcation points in the $(\delta,\,k,\,r)$ space ($n=2$,\,$\beta_0=2.5$). \textbf{Right}- Projection of the surface on the plane $(k,\,\delta)$. Here $\delta$ is denoted by $d$.}}
\end{figure}
\begin{figure}\centering
\includegraphics[width=0.49\linewidth]{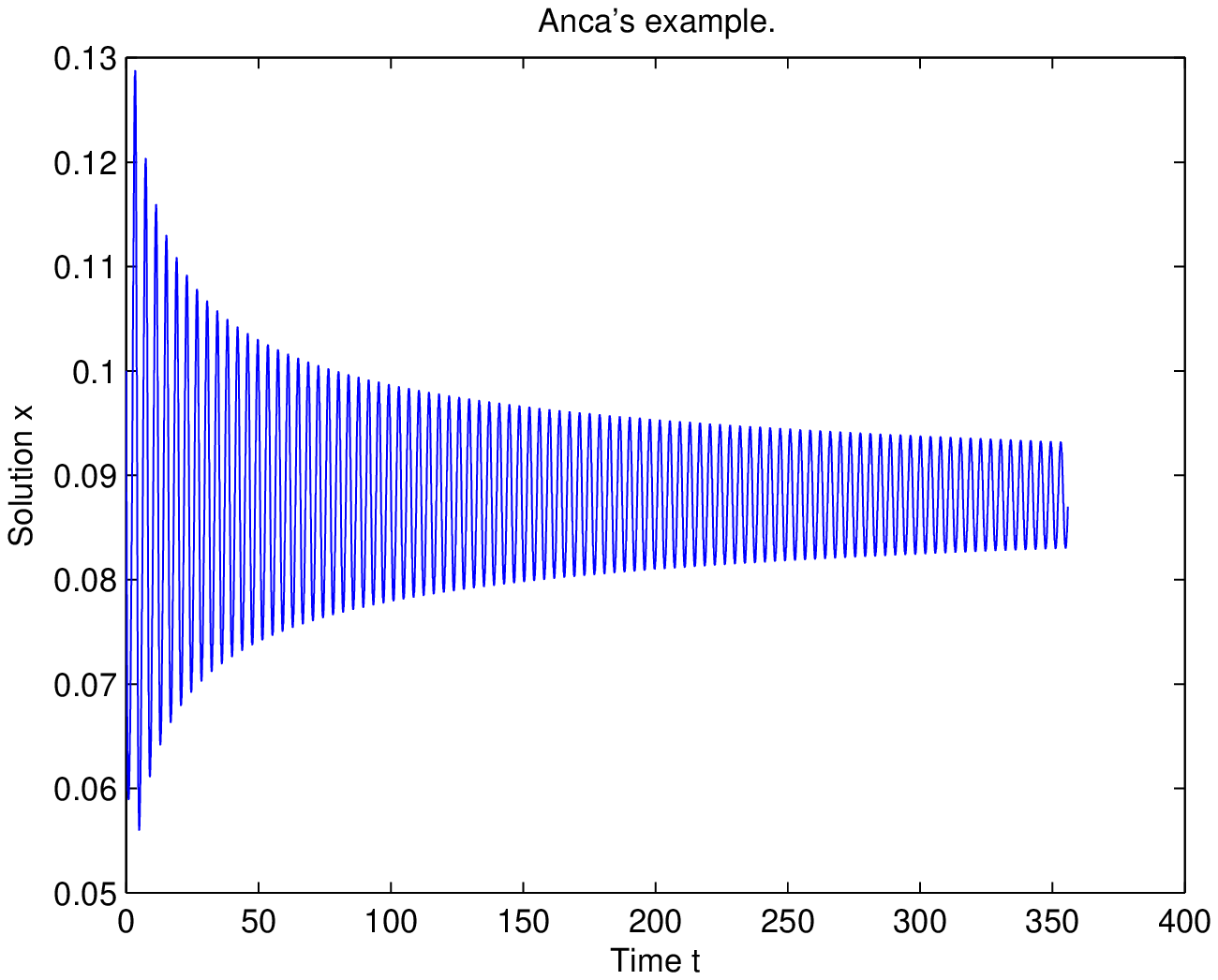}\includegraphics[width=0.49\linewidth]{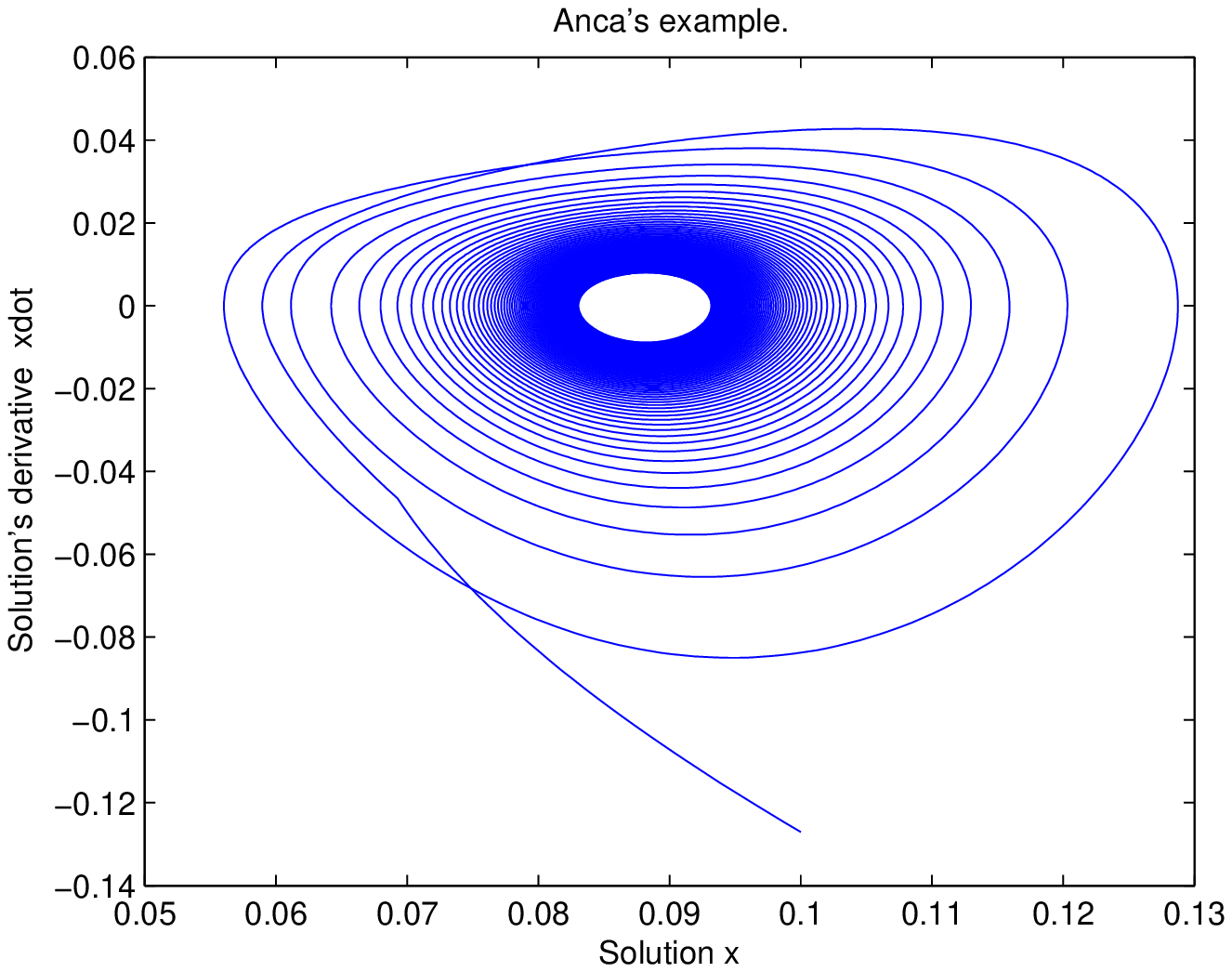}\\\caption{\small{\textbf{Left}- $x(t)$ versus $t$ before the Hopf bifurcation point, at $r=0.3558$. \textbf{Right}- $\dot{x}(t)$ versus $x(t)$ at $r=0.3558$.}}
\includegraphics[width=0.49\linewidth]{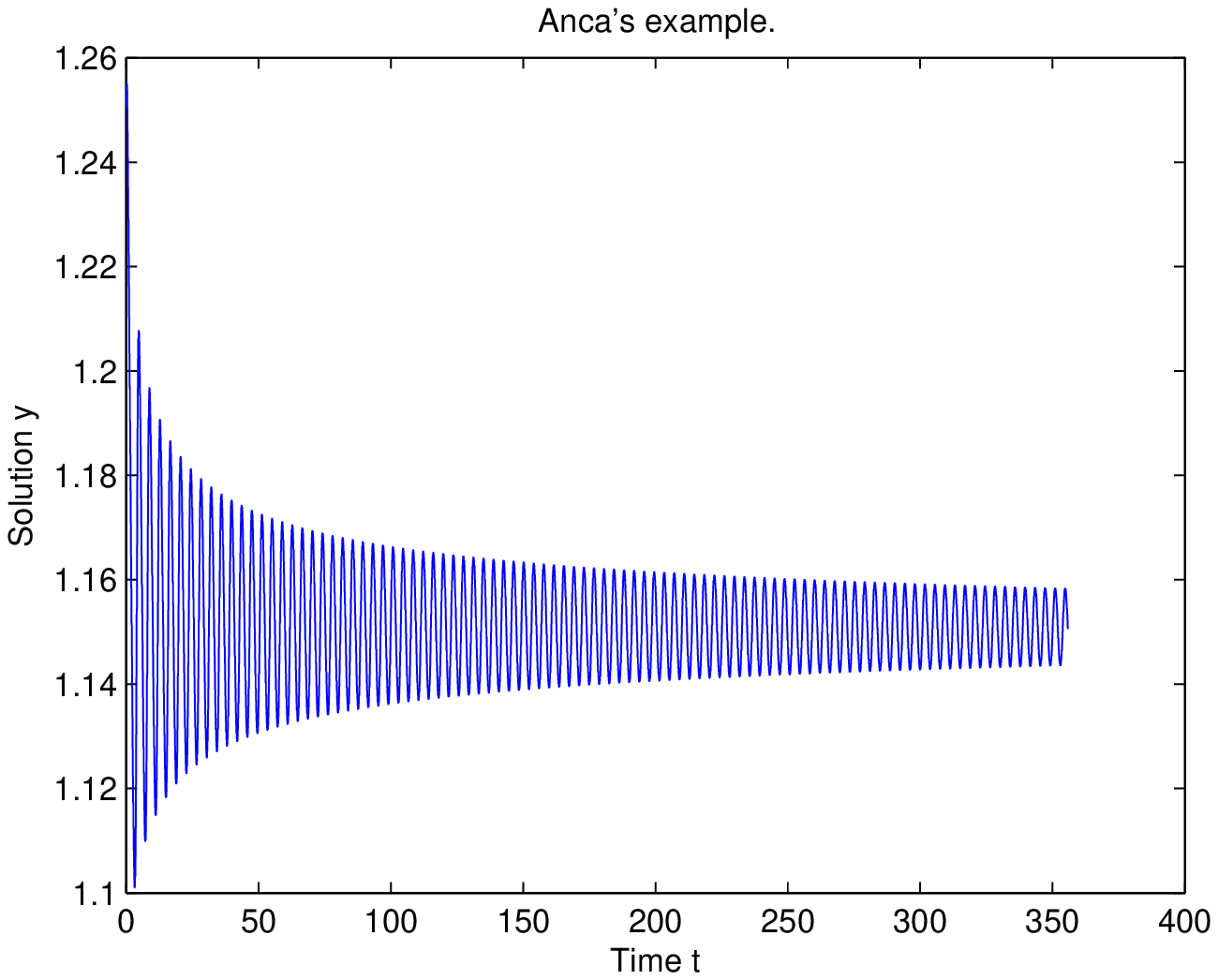}\includegraphics[width=0.49\linewidth]{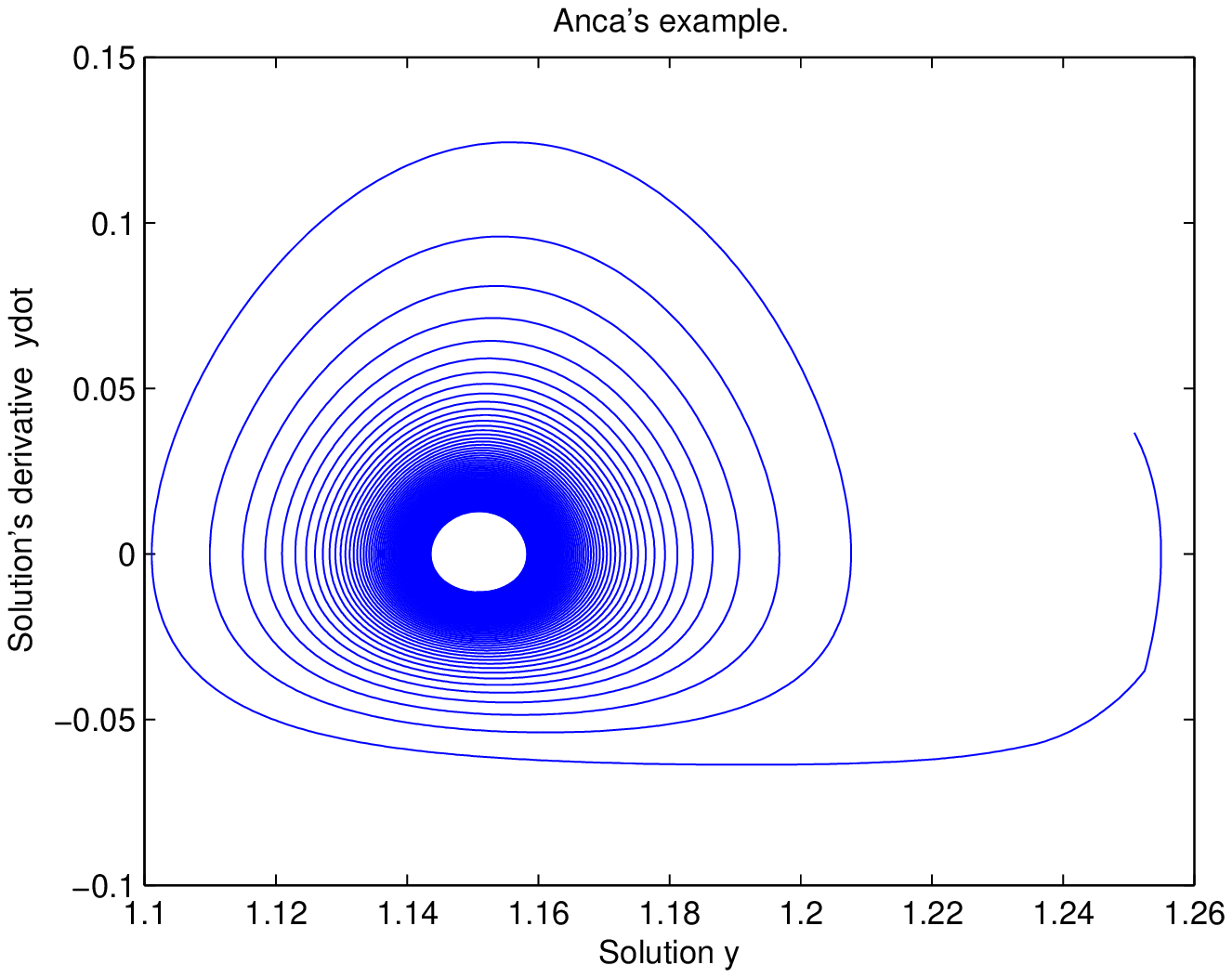}\caption{\small{\textbf{Left}- $y(t)$ versus $t$ before the Hopf bifurcation point, at $r=0.3558$. \textbf{Right}- $\dot{y}(t)$ versus $y(t)$ at $r=0.3558$.}}
\end{figure}

In order to have an image of the points in the parameters space where this happens, we fix two parameters, $n$ and $\beta_0$ and in the three dimensional space of the other parameters, $(k,\,\delta,\,r)$, we represent $r$ as function of $(k,\,\delta)$, by using \eqref{cond-Hopf}. In Fig. 1  such a plot, for $n=1,\,\beta_0=2.5$, can be seen, together with the projection of this surface on the $(k,\,\delta)$ plane (representing the domain of definition of the function in \eqref{cond-Hopf}).

\begin{figure}\centering
\includegraphics[width=0.49\linewidth]{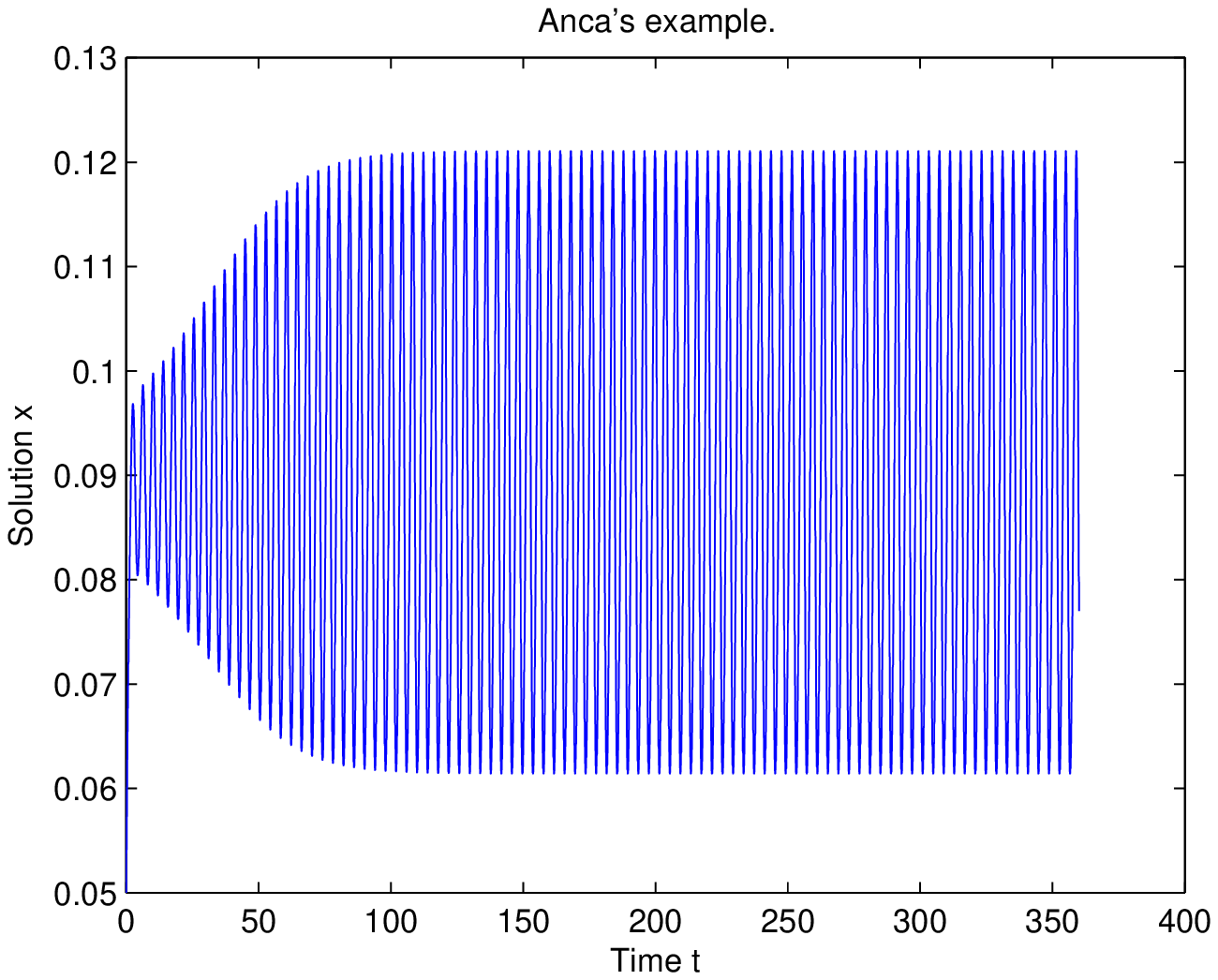}\includegraphics[width=0.49\linewidth]{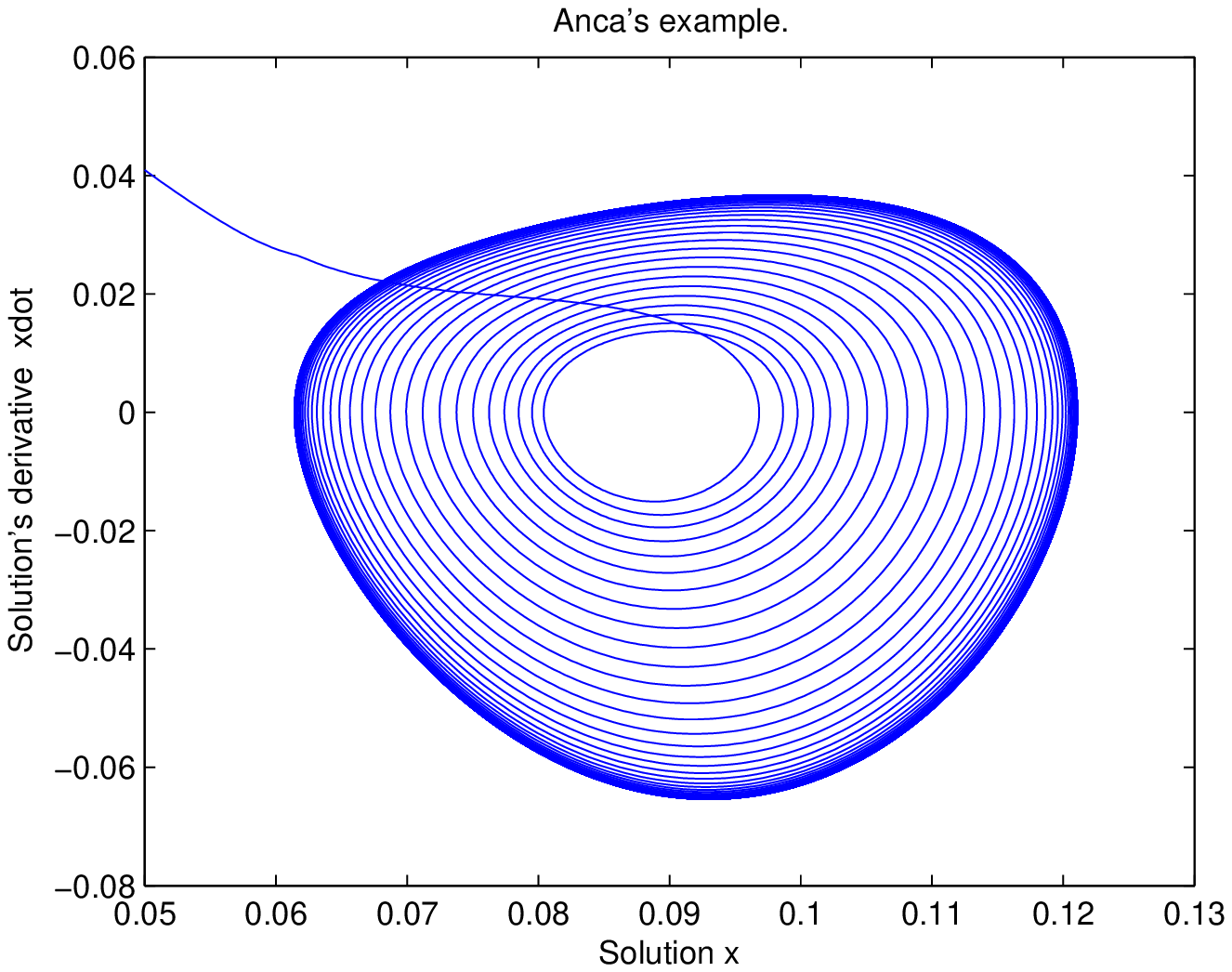}\\\caption{\small{\textbf{Left}- $x(t)$ versus $t$ after the Hopf bifurcation point, at $r=0.36$. \textbf{Right}- $\dot{x}(t)$ versus $x(t)$ at $r=0.36$.}}
\includegraphics[width=0.49\linewidth]{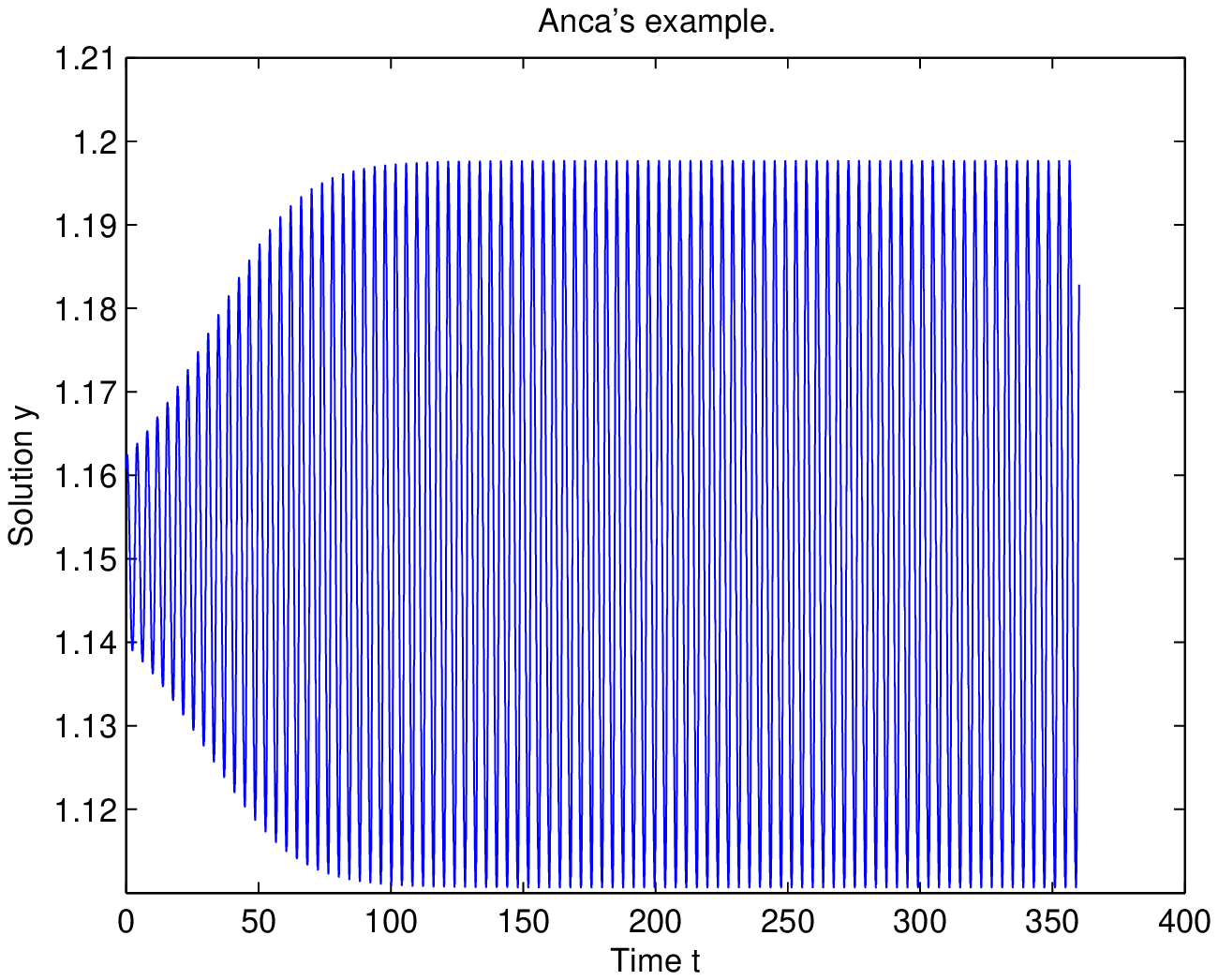}\includegraphics[width=0.49\linewidth]{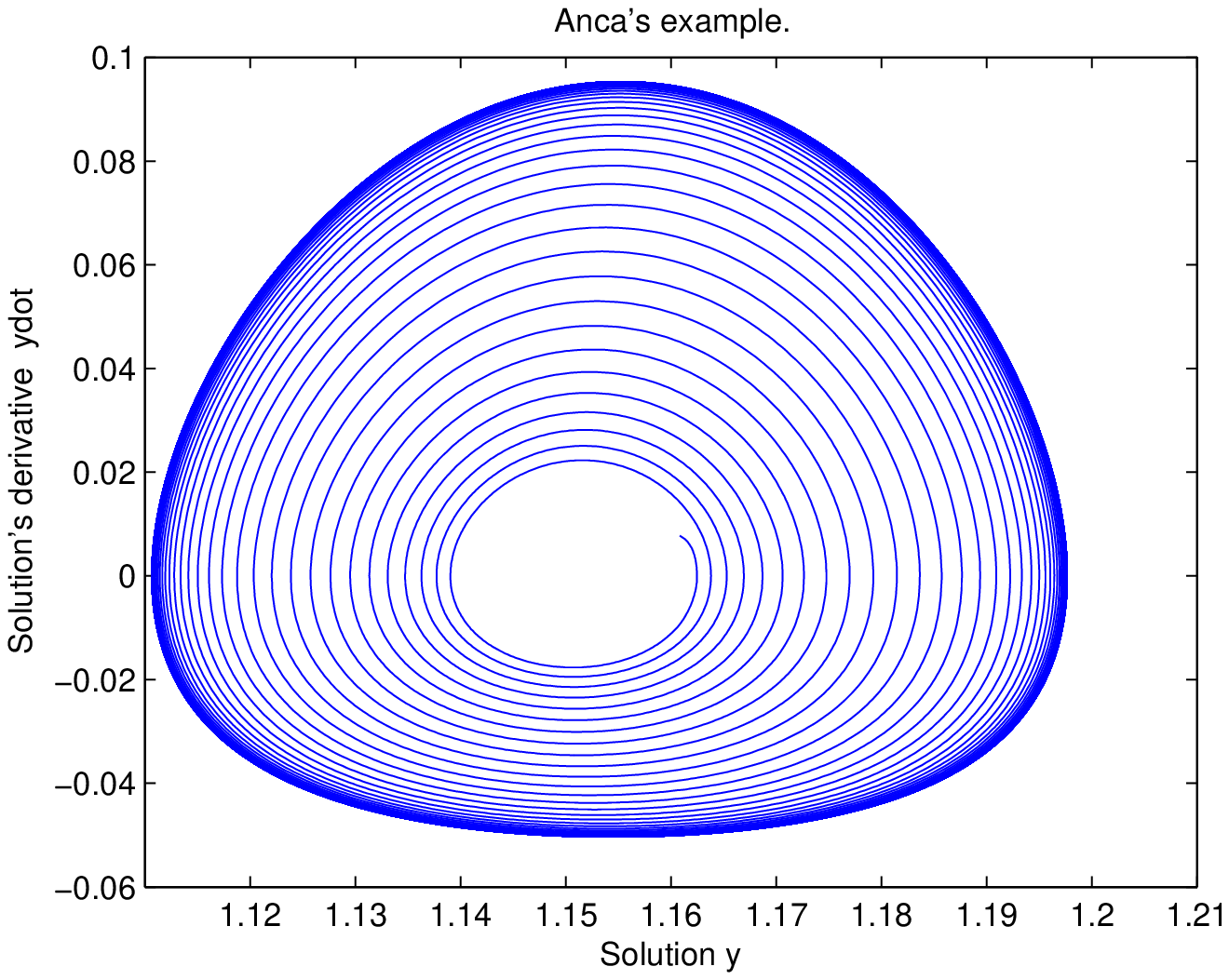}\caption{\small{\textbf{Left}- $y(t)$ versus $t$ after the Hopf bifurcation point, at $r=0.36$. \textbf{Right}- $\dot{y}(t)$ versus $y(t)$ at $r=0.36$.}}
\end{figure}

As a matter of fact, condition \eqref{cond-Hopf} is only a necessary condition for Hopf bifurcation and it must be completed with the condition that the first Lyapunov coefficient, $l_1$, is different of zero (for the definition of $l_1$ see \cite{AI-Hopf}, where the restriction of the problem to a center manifold is used). If $l_1=0$, then depending on the number of varying parameters, either a Hopf degenerate bifurcation (a single varying parameter) or a Bautin bifurcation (two varying parameters) takes place.

\vspace{0.3cm}

In Figures 2, 3, 4, 5, we present the solutions  $x(t)$ and $y(t)$
at a typical Hopf bifurcation point,
($n=12,\,\beta_0=1.77,\,\delta=0.05,\,k=1.18074,\,r=0.3559114$).
Here $r$ was taken the bifurcation parameter. Figures 2 and 3
present the solutions, for the parameter $r$ at the left of the
bifurcation value (at $r=0.3558$) while Figures 4 and 5 present the
solutions, for $r$ at the right of the bifurcation value ($r=0.36$).
It is clear that the behaviour of $y$ determines that of $x$. We
remark that the first Lyapunov coefficient for this Bautin
bifurcation is $l_1=-43.71.$

\section{Points of Bautin bifurcation for $y$. \\
Influence on $x$}

In \cite{AVI-RMG} we studied the Bautin bifurcation for eq.
\eqref{eq-y}. The Bautin bifurcation  occurs when $l_1=0$ and the
second Lyapunov coefficient, $l_2$, is different of zero. For the
definition of $l_2$ see \cite{AVI-RMG} - where again the problem is
restricted to a two-dimensional center manifold. By this restriction
a two dimensional system of differential equations is obtained. For
this kind of problems, the Bautin bifurcation was studied in
\cite{Ku}.

In the Bautin bifurcation two parameters are varied, hence we may
speak of a \textit{parameters plane} and the normal form of this
bifurcation in polar coordinates is
\[\dot{\rho}=\rho(b_1+b_2\rho^2+\zeta\rho^4),
\]
\[\dot{\theta}=1,
\]
where $\zeta=\mathrm{sgn} (l_2).$

The Bautin bifurcation is characterized by the fact that, if we vary
the parameters in the parameters plane, around the bifurcation
point, we can meet four possible situations.  These are, as can be
seen on the bifurcation diagram, Fig. 6: an equilibrium point  (in
the domain denoted by 1 of the parameters plane), an equilibrium
surrounded  by a limit cycle (in domain 2),  an equilibrium
surrounded by two cycles one inside the other (zone 3), and,
finally, an equilibrium that is attracting on one side and repelling
on its other side (for points on the curve T). The curve T is the
curve on which the discriminant of equation $b_1+b_2\eta+\zeta\eta^2=0$
is equal to zero and the solution $\eta$ is positive (since $\eta=\rho^2$).

The stability properties of the trajectories described above depend
on the sign of $l_2$, as we can see in the bifurcation diagrams
(Fig. 6).

The letter ``H'' on the bifurcation diagram shows that, while crossing
the vertical axis, in the sense indicated by the small arrows, Hopf bifurcation takes place.

\vspace{0.3cm}

\begin{figure}\centering
\includegraphics[width=0.45\linewidth]{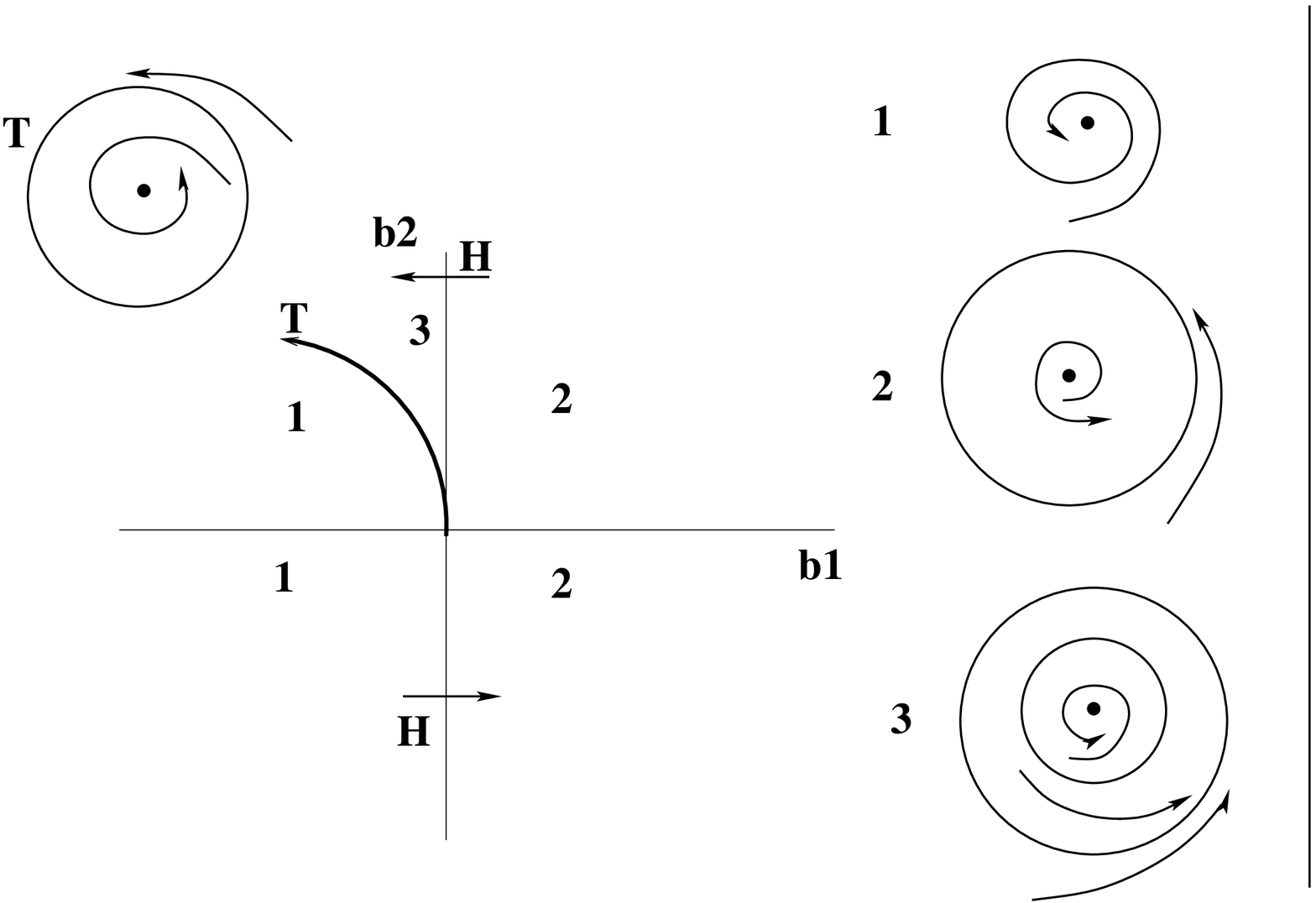}\includegraphics[width=0.42\linewidth]{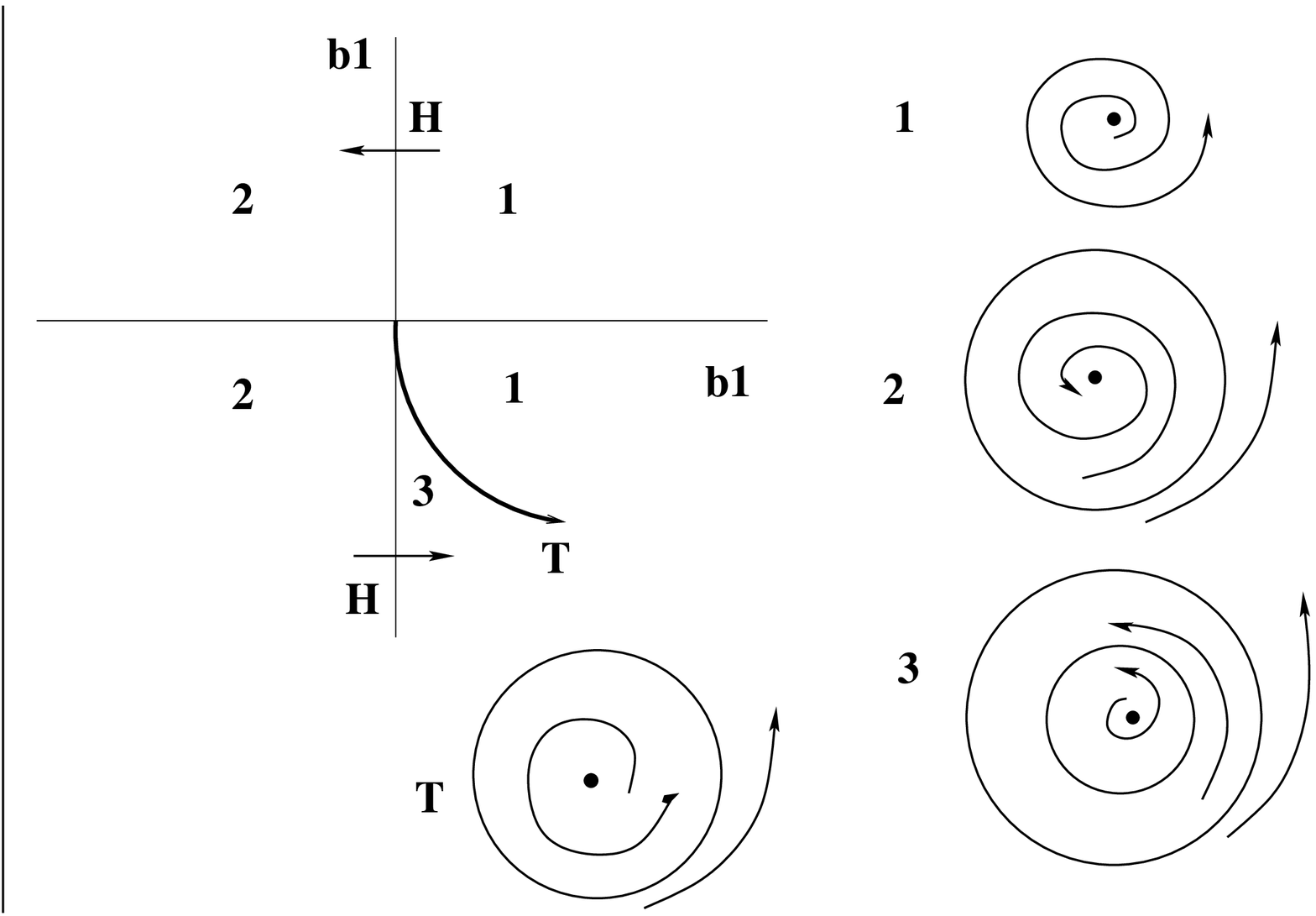}\caption{Bautin bifurcation diagram. \textbf{Left}- the case $l_2<0$; \textbf{Right}- the case $l_2>0$.}
\end{figure}

The algorithm for finding Bautin bifurcation points for eq. \eqref{eq-y} was described in detail in \cite{AVI-RMG}. We only shortly remind the results of our search.

\vspace{0.3cm}

We considered the zone of parameters $(n,\,\beta_0,\,k)$ given by:
$$(n,\,\beta_0,\,k)\in\{1,1.5,2,3,...,12\}\times\{0.5,1,1.5,2,2.5\}\times \{1.1,1.2,...,1.9\}
$$ (we remind that $k>1$ in order that $x_2$ exist, and $k<2$). $(n,\beta_0,k)$ were
taken as fixed parameters, and $(\delta,\,r)$ as variable
parameters. Hence, for every $(n^*,\beta^*_0,k^*)$ fixed, we looked
for a $\delta^*$ and a $r^*$  such that $Re
\lambda_{1,2}(\alpha^*)=0$ and $l_1(\alpha^*)=0,$ where we denoted
by $\alpha$ the vector of parameters, \break
$\alpha=(n,\,\beta_0,\,k,\,\delta,\,r)$.

\vspace{0.3cm}

For $\mathbf{n=1}$ we showed \cite{AVI-RMG} that the condition \eqref{cond-Hopf} for Hopf bifurcation is not satisfied, hence no Bautin bifurcation can occur.

\vspace{0.3cm}

For  $\mathbf{n=1.5}$ and $\mathbf{n=2},$ and for each
$\beta_0\in\{0.5,\,1,\,1.5,\,2,\,2.5 \},$ and each
$k\in\{1.1,\,1.2,...,1.9\}$ we found Hopf points with $l_1=0$. For
these points we computed the second Lyapunov coefficient and found
that for each of the bifurcation points obtained, $l_2<0$. Hence
they are Bautin bifurcation points and the bifurcation diagram is
that of Fig. 6, left.

\vspace{0.3cm}

For $\mathbf{3\leq n \leq 12}$ and
$\beta_0\in\{0.5,\,1,\,1.5,\,2,\,2.5 \},$
$k\in\{1.1,\,1.2,...,1.9\}$, the first Lyapunov coefficient was find
to be always negative, hence, no Bautin bifurcation is possible for
this range of parameters.

\vspace{0.3cm}

We give,  below, the points, found by us, in the parameters space
where $l_1=0$ for the case $n=2$ and $\beta_0\in \{0.5,\,1,\,1.5,\,2
\}$ and $k$ as above. Also we present the plots of these Bautin
bifurcation points on the surfaces of Hopf points (Figs. 7 - 10). By
connecting these points, we obtain approximate curves of Bautin
bifurcation points on the surfaces of Hopf points.

A table and a figure of the same type, for the case
$n=2,\,\beta_0=2.5$, were published in \cite{AVI-RMG}.

\begin{figure}[h]\centering
\begin{minipage}[h]{0.49\linewidth}
{\tiny\bf\begin{tabular}{|c|c|c|c|} \hline\
 k  & $\delta$ & $r$ & $l_2$    \\ \hline
 1.1  & 0.0045705962 & 26.125314 & -0.021 \\ \hline
 1.2  & 0.0090491351 & 25.751524 & -0.0151    \\ \hline
 1.3  & 0.0134437887 & 25.422162 & -0.0124    \\ \hline
 1.4  & 0.0177612407 & 25.130258 & -0.0108   \\ \hline
 1.5  & 0.0220070315 & 24.870352 & -0.0097   \\ \hline
 1.6  & 0.0261858065 & 24.638093 & -0.0088    \\ \hline
 1.7  & 0.0303014988 & 24.429962 & -0.0081    \\ \hline
 1.8  & 0.0343574676 & 24.243076 & -0.0076    \\ \hline
 1.9  & 0.0383566021 & 24.075039 & -0.0071    \\ \hline
\end{tabular}}
\end{minipage}
\begin{minipage}[h]{0.49\linewidth}
\includegraphics[width=\linewidth]{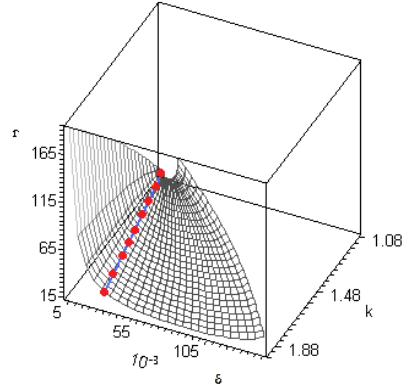}
\end{minipage}
\caption{Hopf codimension two points for
$n=2,\,\beta_0=0.5.$}\end{figure}

\begin{figure}[h]\centering
\begin{minipage}[h]{0.49\linewidth}
{\tiny\bf\begin{tabular}{|c|c|c|c|} \hline\
       k  & $\delta$ & $r$ & $l_2$    \\ \hline
 1.1 & 0.0091411924 &13.062657 & -0.0205 \\ \hline
 1.2 & 0.0180982702  &12.875762 & -0.0142\\ \hline
 1.3 & 0.0268875774 &12.711081 & -0.0114   \\ \hline
 1.4 & 0.0355224814  &12.565129 & -0.0097   \\ \hline
 1.5 & 0.0440140630  &12.435176 & -0.0085   \\ \hline
 1.6 & 0.0523716129  &12.319046 & -0.0076  \\ \hline
 1.7 & 0.0606029975  &12.214981 & -0.0069   \\ \hline
 1.8 & 0.0687149345  &12.121538 & -0.0063   \\ \hline
 1.9 & 0.0767132043  &12.037519 &  -0.0059   \\ \hline
\end{tabular}}
\end{minipage}
\begin{minipage}[h]{0.49\linewidth}\includegraphics[width=\linewidth]{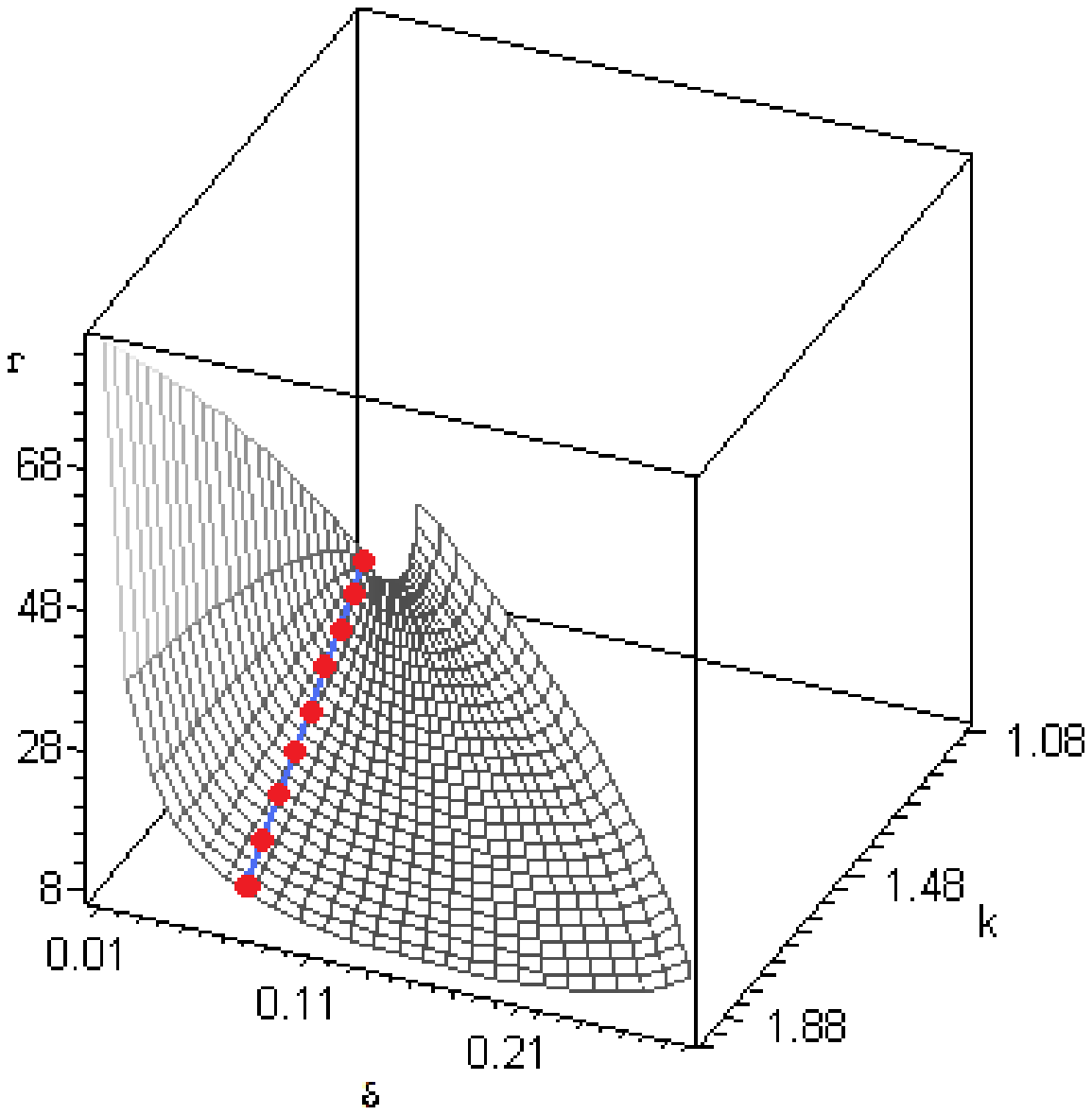}
\end{minipage}
\caption{Hopf codimension two points for
$n=2,\,\beta_0=1.$}\end{figure}

\begin{figure}[h]\centering
\begin{minipage}[h]{0.49\linewidth}
{\tiny\bf\begin{tabular}{|c|c|c|c|} \hline\
 k  & $\delta$ & $r$ & $l_2$    \\ \hline
 1.1  & 0.0137117887 & 8.708438 &  -0.0204   \\ \hline
 1.2  & 0.0271474053 & 8.583841 & -0.0140   \\ \hline
 1.3  & 0.0403313662 & 8.474054 & -0.0112   \\ \hline
 1.4  & 0.0532837222 & 8.376752 & -0.0095   \\ \hline
 1.5  & 0.0660210946 & 8.290117 &  -0.0083   \\ \hline
 1.6  & 0.0785741932 & 8.212697 & -0.0074    \\ \hline
 1.7  & 0.0909044966 & 8.143320 &  -0.0067   \\ \hline
 1.8  & 0.1030724022 & 8.081025 & -0.0061    \\ \hline
 1.9  & 0.1150698062 & 8.025013 & -0.0056 \\ \hline
\end{tabular}}\end{minipage}
\begin{minipage}[h]{0.49\linewidth}\includegraphics[width=\linewidth]{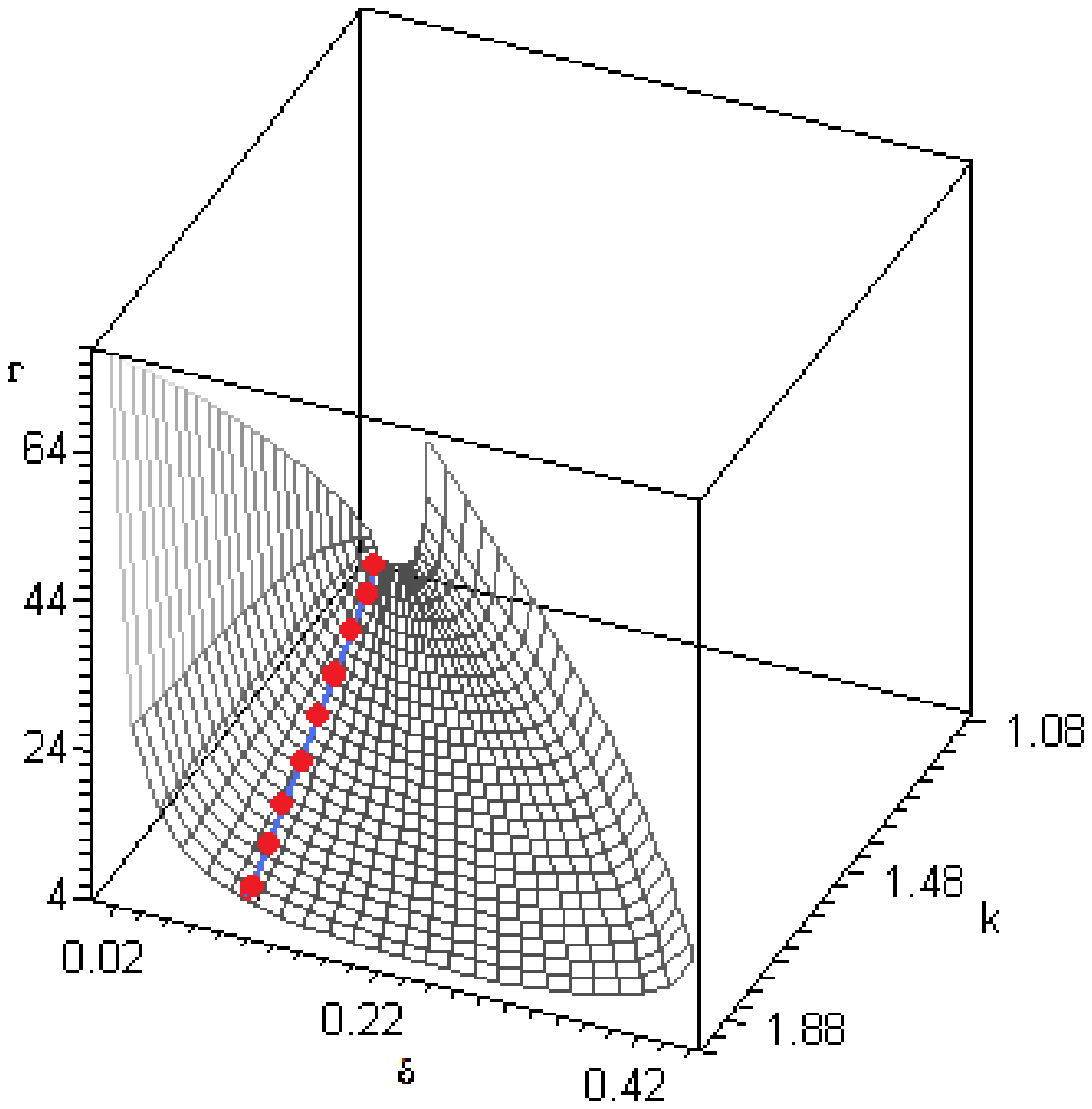}
\end{minipage}
\caption{Hopf codimension two points for
$n=2,\,\beta_0=1.5.$}\end{figure}

\begin{figure}[h]\centering
\begin{minipage}[h]{0.49\linewidth}
{\tiny\bf\begin{tabular}{|c|c|c|c|} \hline\
 k  & $\delta$ & $r$ & $l_2$    \\ \hline
 1.1  & 0.018282385 & 6.531328 & -0.0203   \\ \hline
 1.2  & 0.036196540 & 6.437880 & -0.014  \\ \hline
 1.3  & 0.053775154 & 6.355540 & -0.0111   \\ \hline
 1.4  & 0.071044963 & 6.282564 & -0.0093  \\ \hline
 1.5  & 0.088028126 & 6.217588 & -0.0082   \\ \hline
 1.6 & 0.104743225 & 6.159523 & -0.0073   \\ \hline
 1.7  & 0.121205995 & 6.107490 & -0.0066    \\ \hline
 1.8  & 0.137429869 & 6.060769 &  -0.0060   \\ \hline
 1.9  & 0.153426408 & 6.018759 & -0.0055    \\ \hline
\end{tabular}}\end{minipage}
\begin{minipage}[h]{0.49\linewidth}\includegraphics[width=\linewidth]{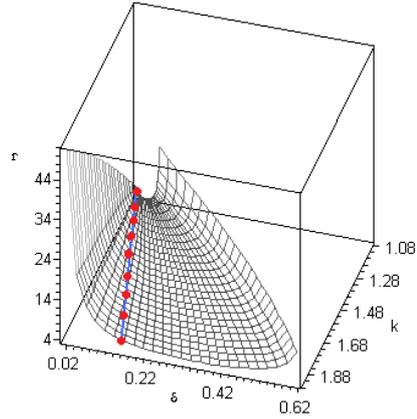}
\end{minipage}
\caption{Hopf codimension two points for
$n=2,\,\beta_0=2.$}\end{figure}

\vspace{0.3cm}

In \cite{AVI-RMG}, in order to confirm by numerical integration the
behaviour predicted by the bifurcation diagram, we considered a
Bautin bifurcation point $P^*$, identified by the algorithm in
\cite{AVI-RMG}, at $n^*=2,\,\beta^*_0=2.5,$ $k^*=1.01,\,$ and
$r^*=5.301432998$, $\delta^*=0.0023073665$.

\vspace{0.3cm}

\begin{figure}\centering
\includegraphics[width=0.49\linewidth]{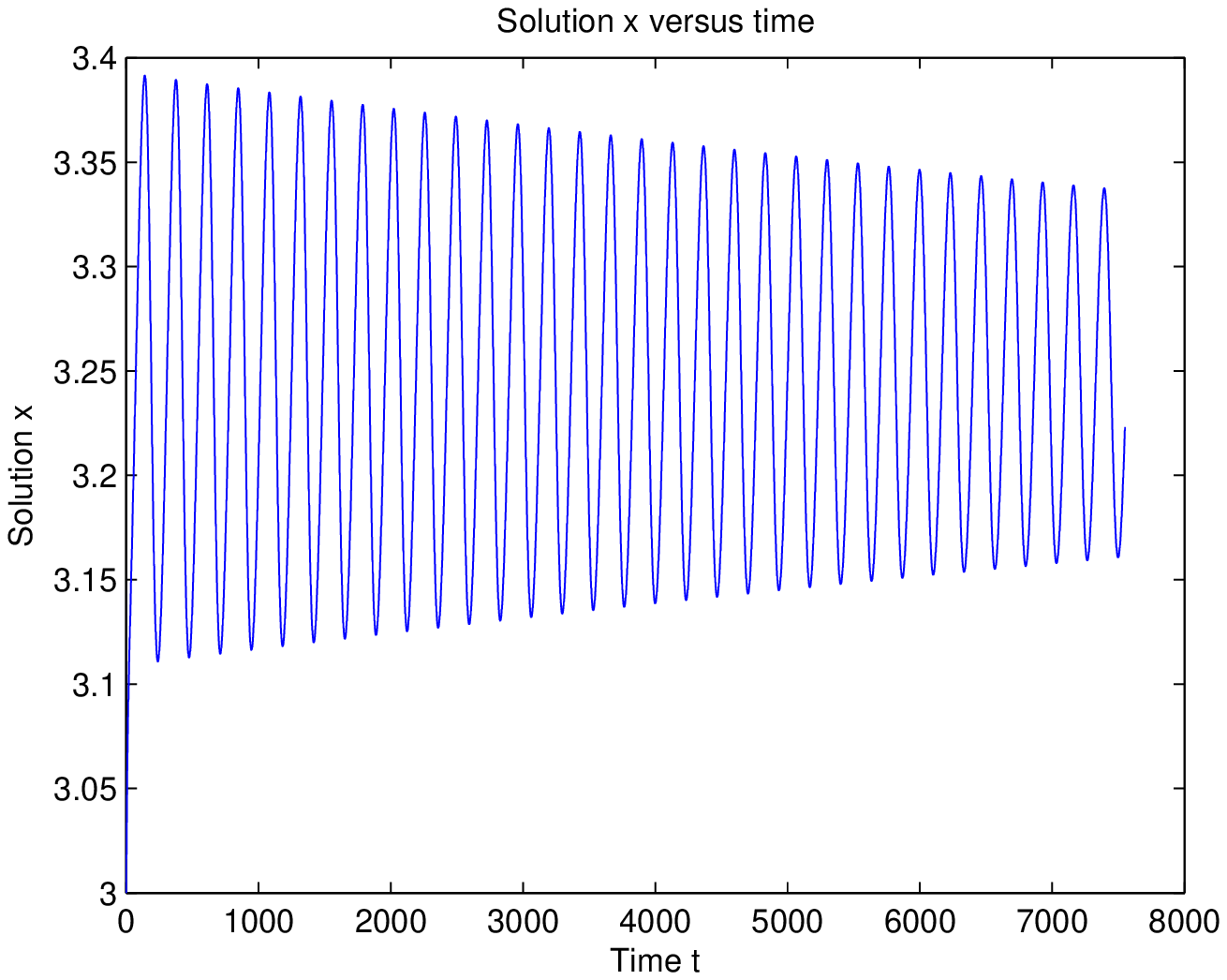}\includegraphics[width=0.49\linewidth]{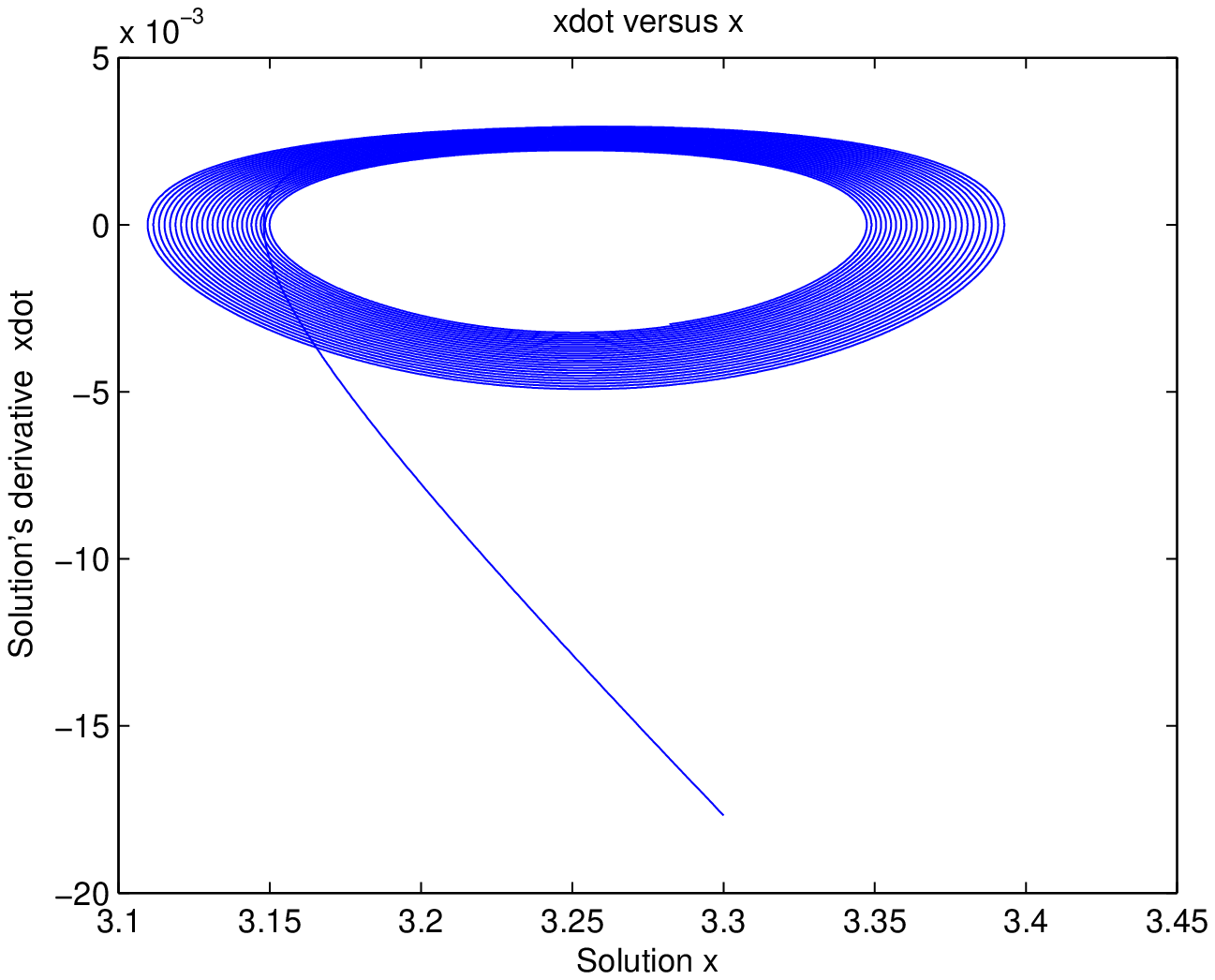}\\
\includegraphics[width=0.49\linewidth]{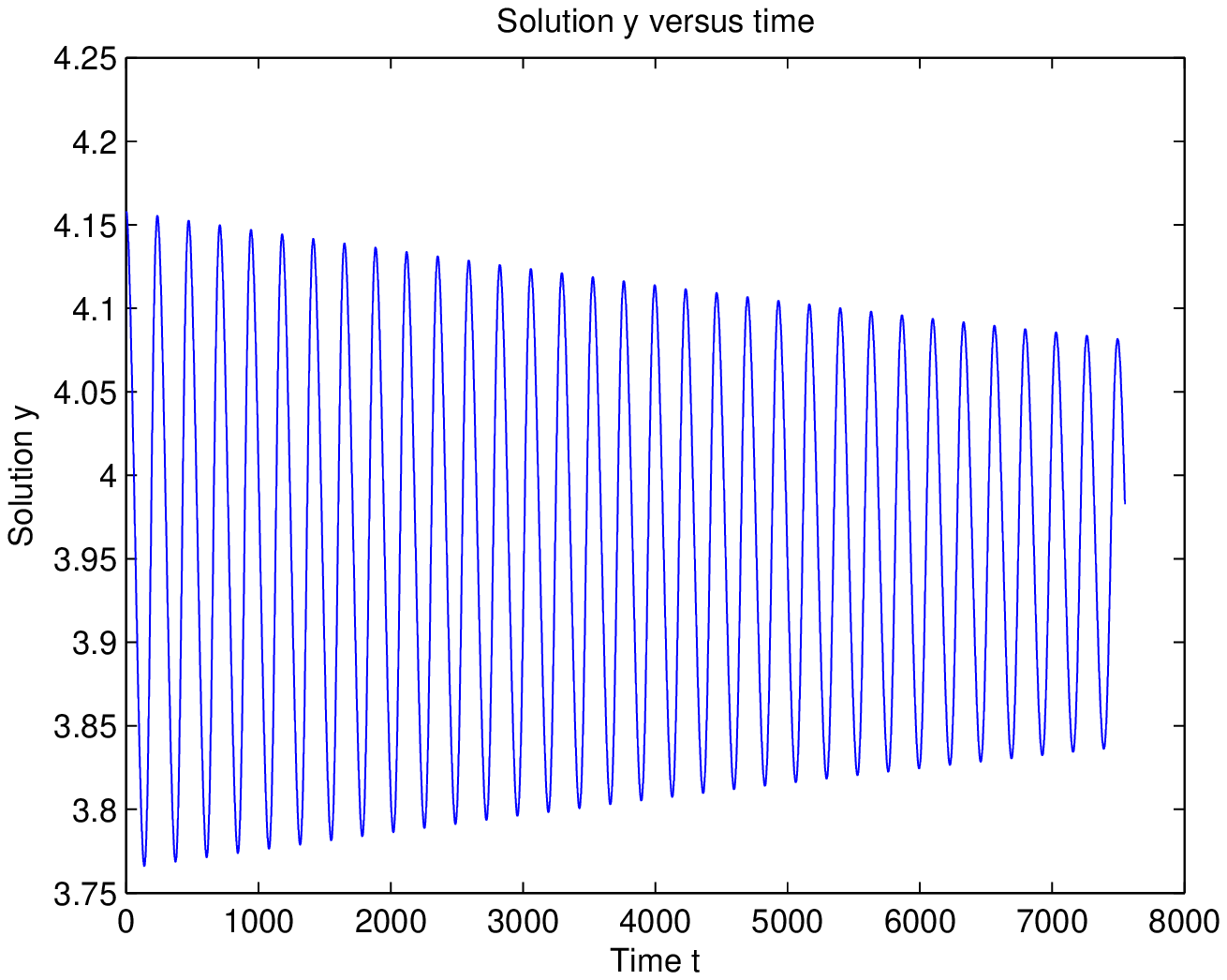}\includegraphics[width=0.49\linewidth]{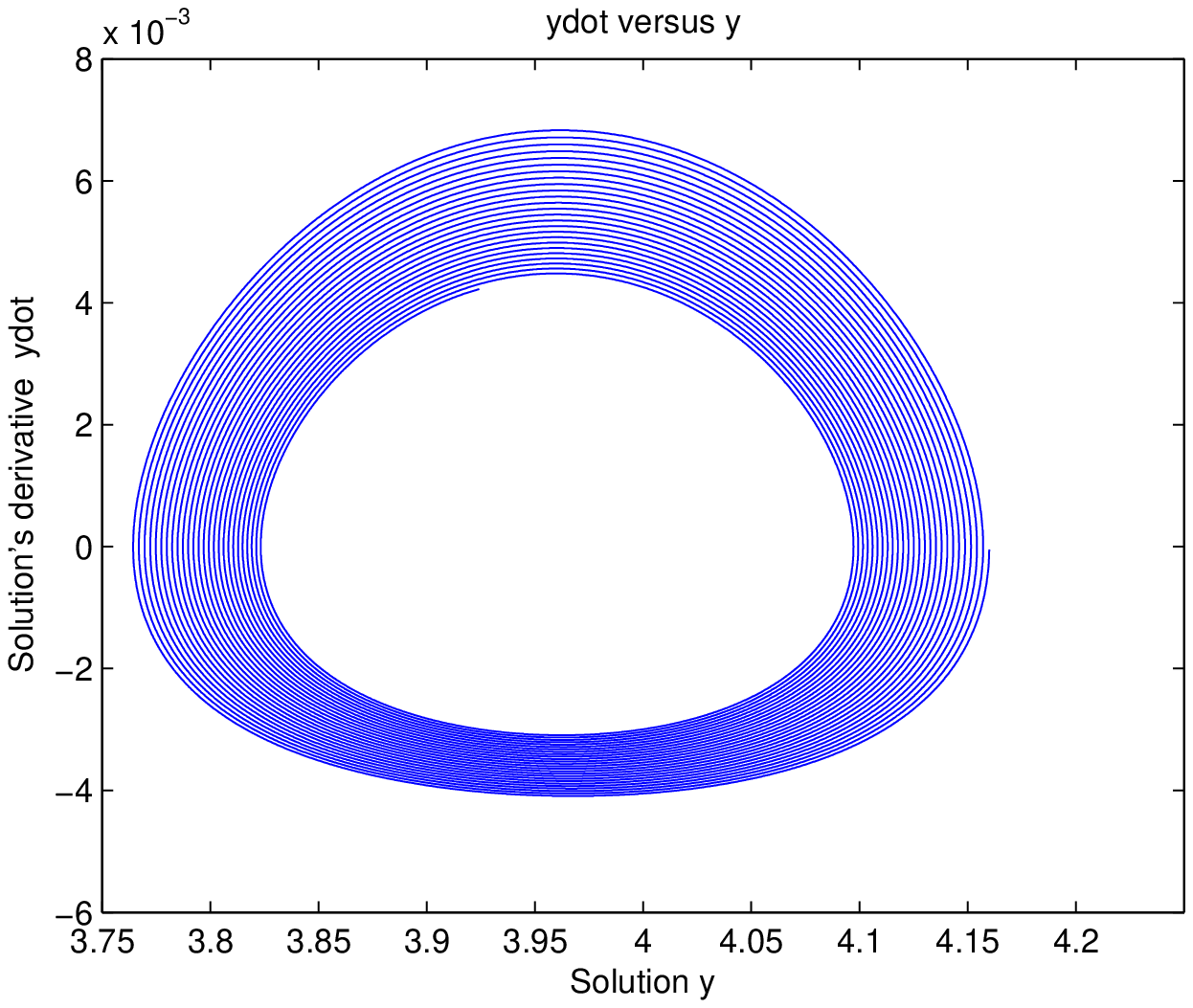}
\caption{Graphs of $x(\cdot),$ respectively $y(\cdot)$, for point $P_3,\;\;c=0.2$ .}
\end{figure}

\begin{figure}\centering
\includegraphics[width=0.49\linewidth]{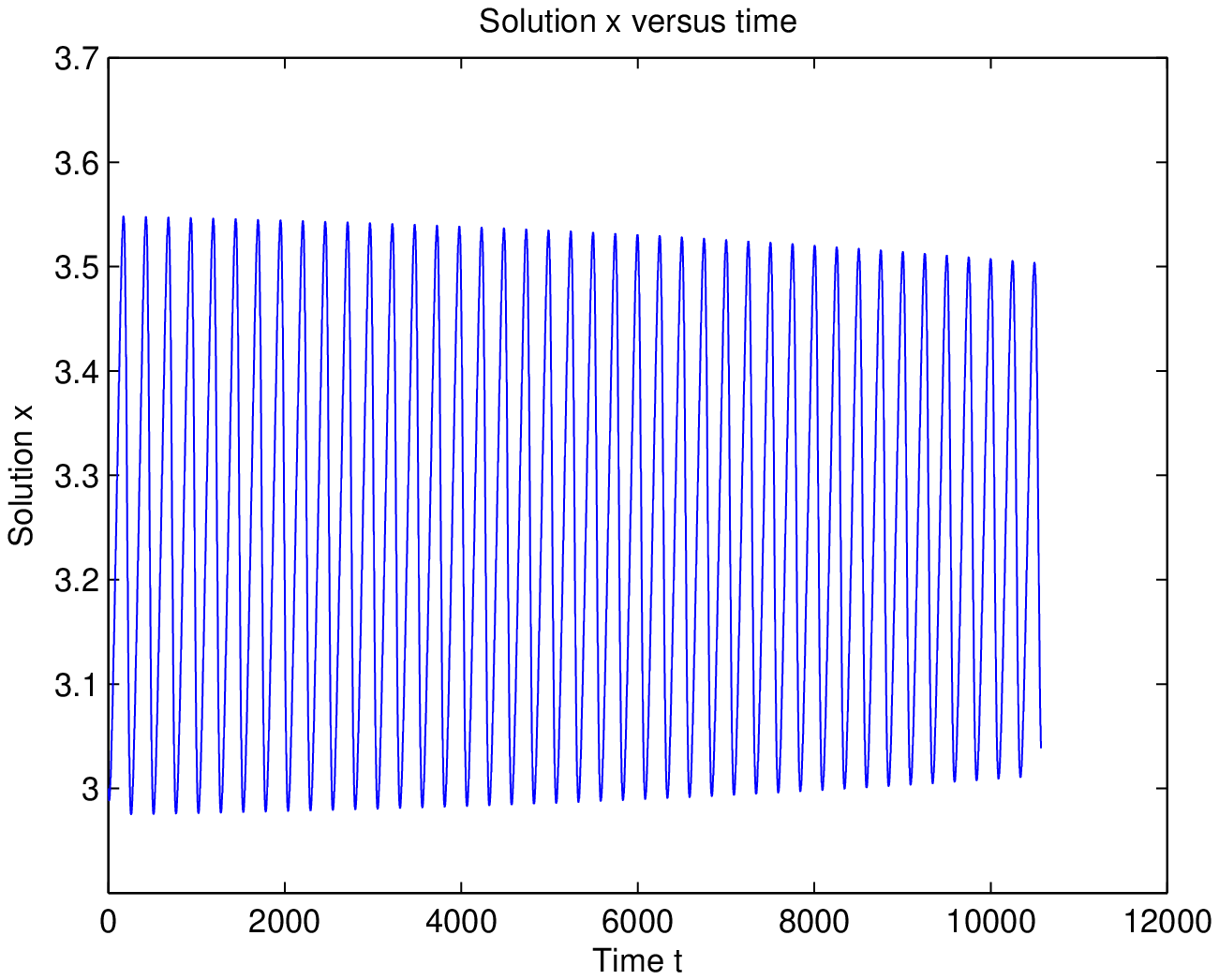}\includegraphics[width=0.49\linewidth]{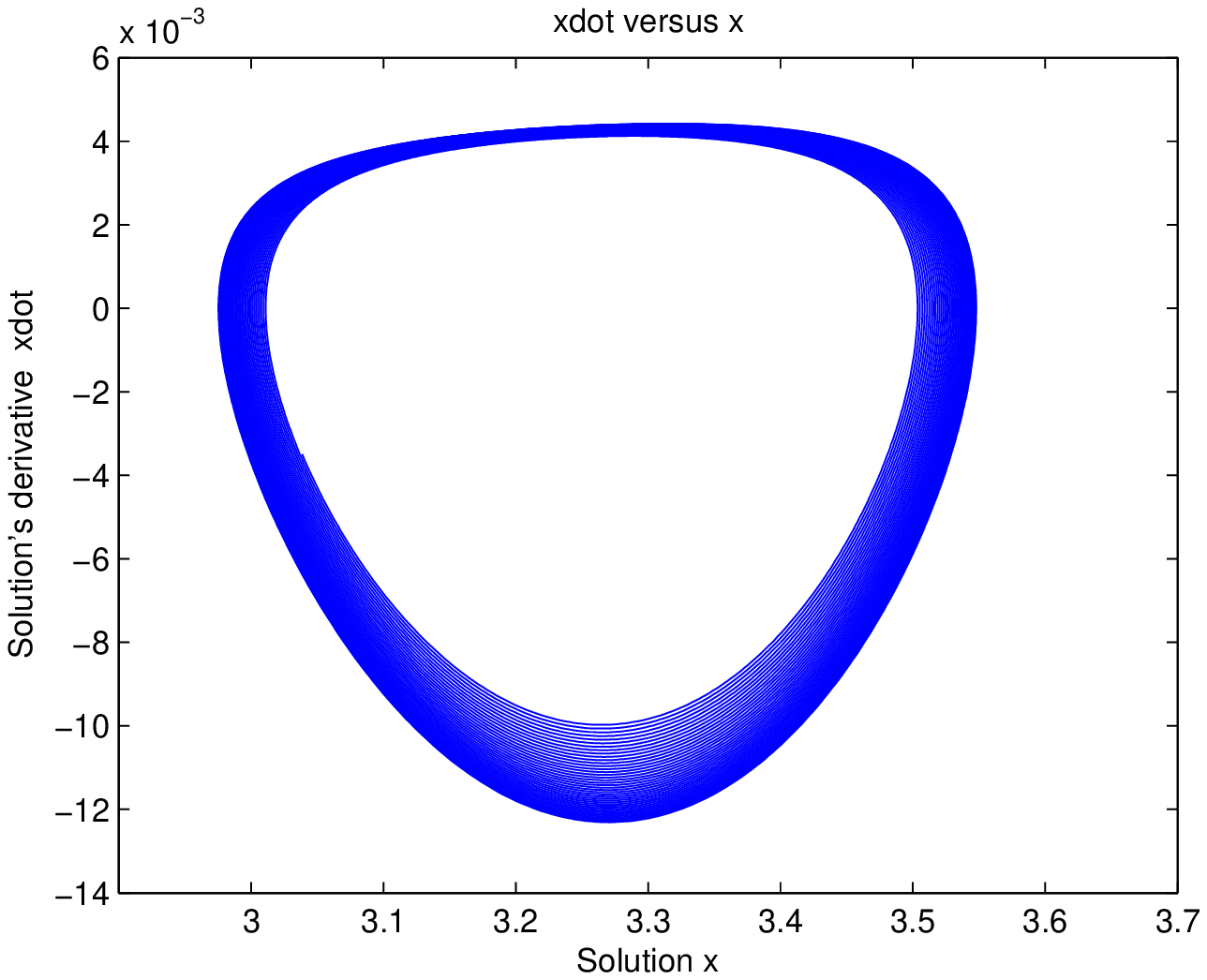}\\
\includegraphics[width=0.49\linewidth]{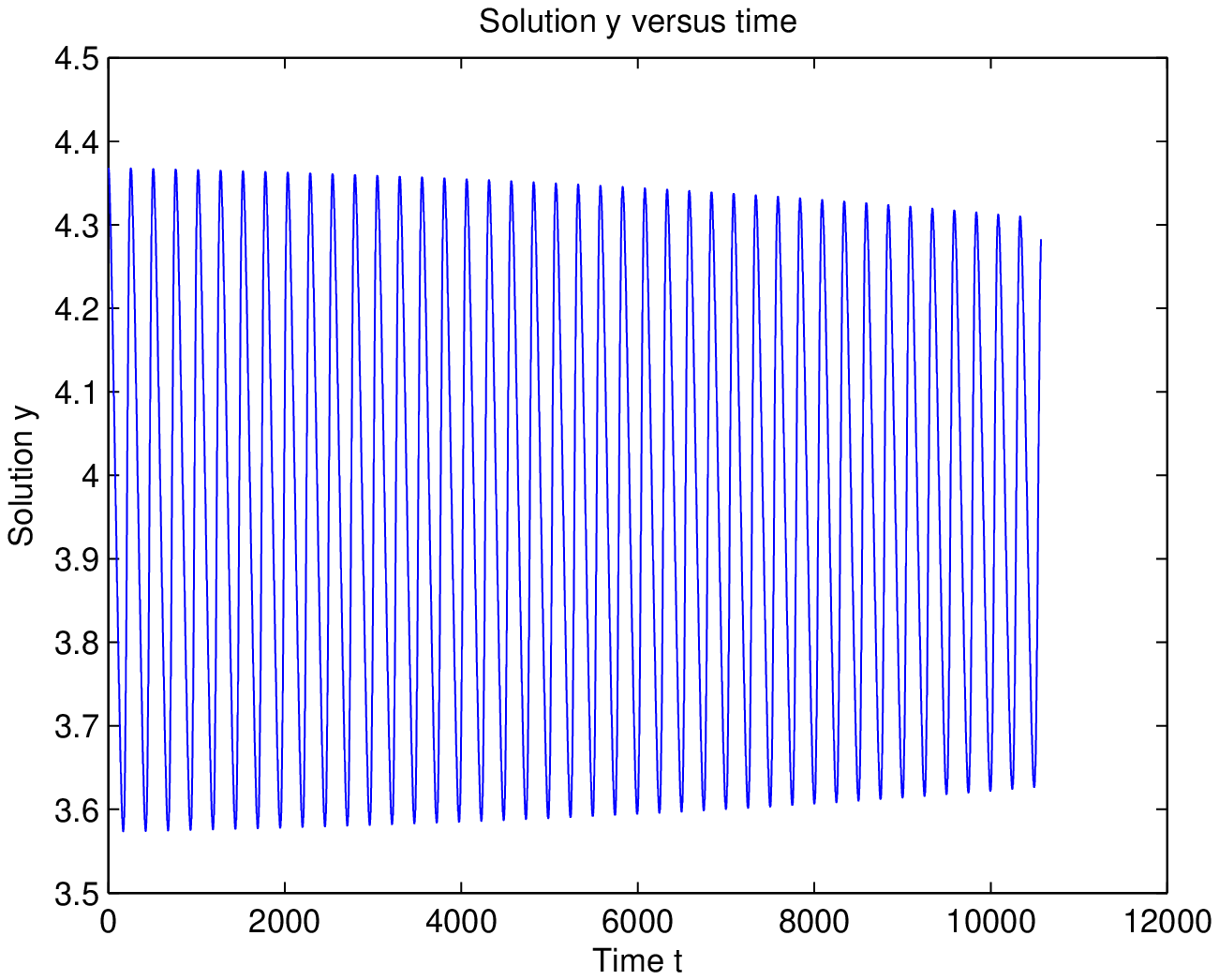}\includegraphics[width=0.49\linewidth]{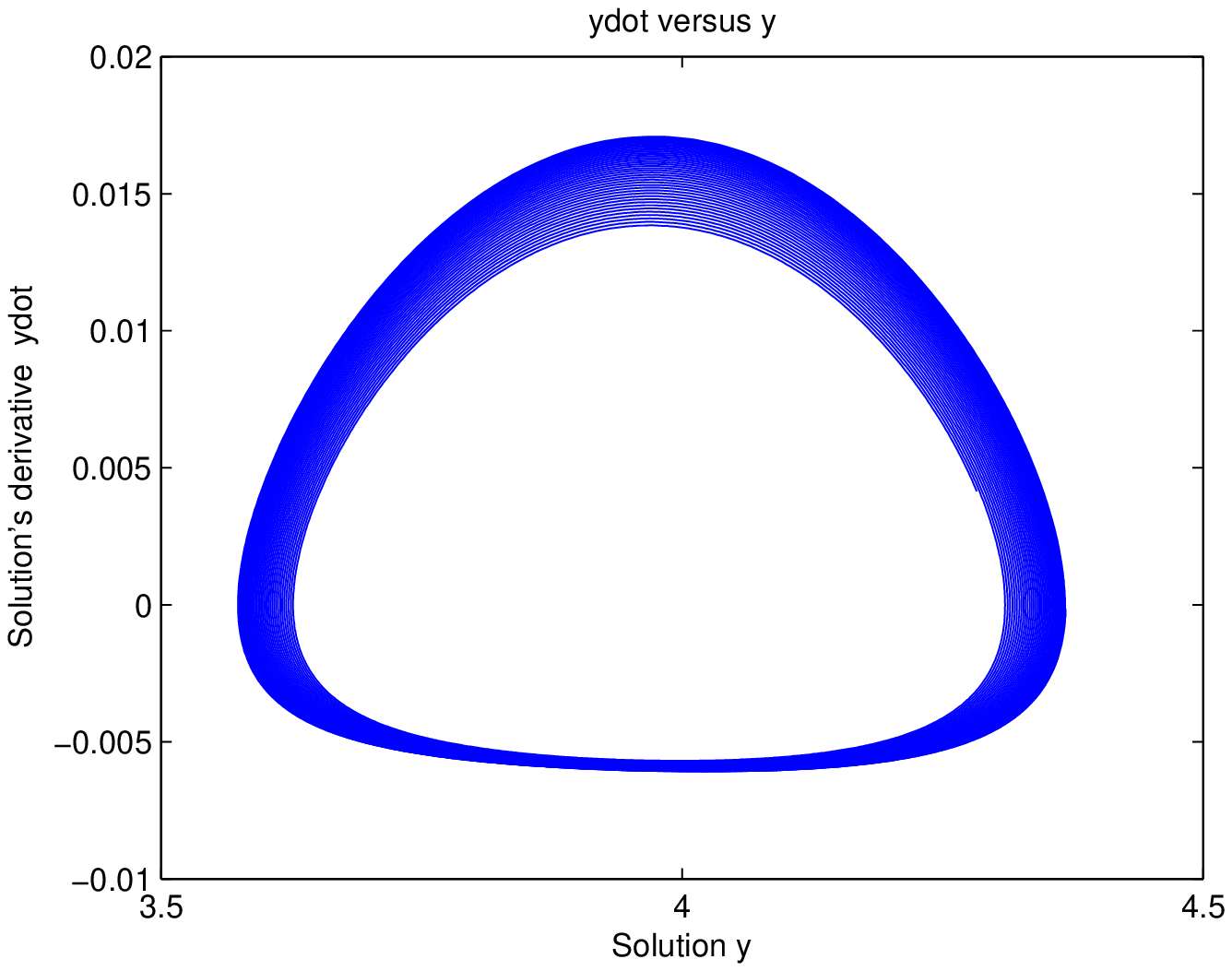}
\caption{Graphs of $x(\cdot),$ respectively $y(\cdot)$, for point
$P_3,\;\;c=0.41$ .}
\end{figure}

\begin{figure}\centering
\includegraphics[width=0.49\linewidth]{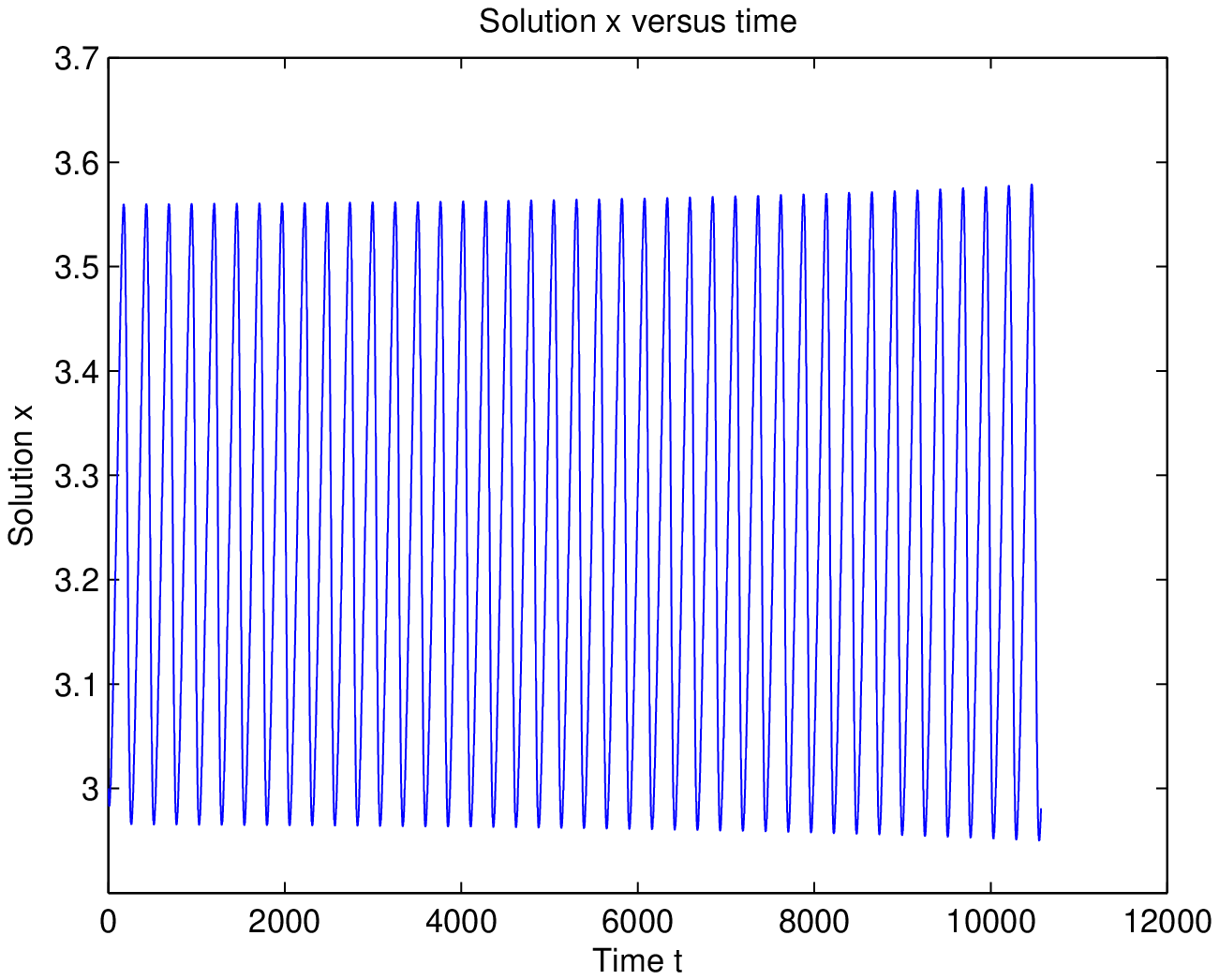}\includegraphics[width=0.49\linewidth]{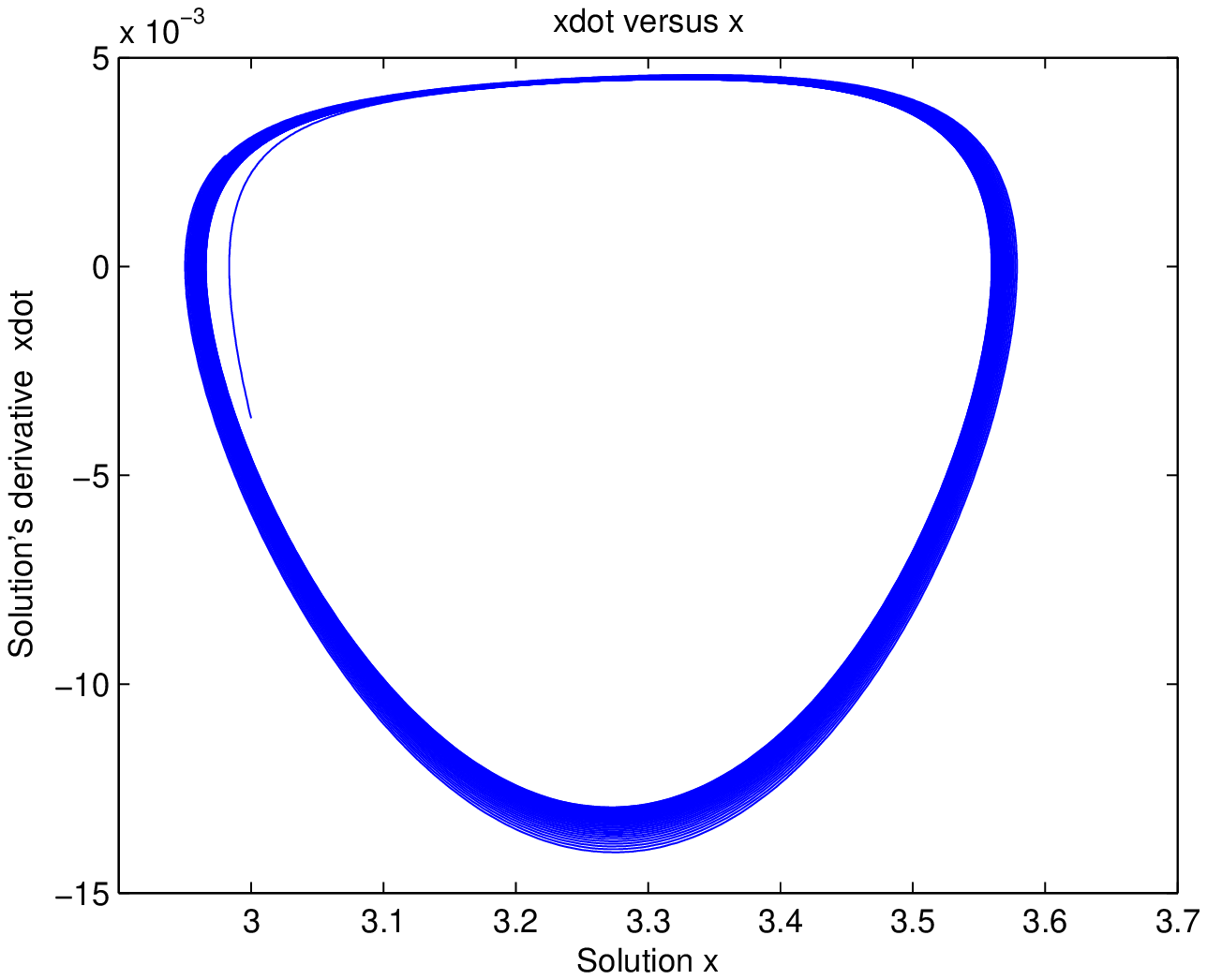}\\
\includegraphics[width=0.49\linewidth]{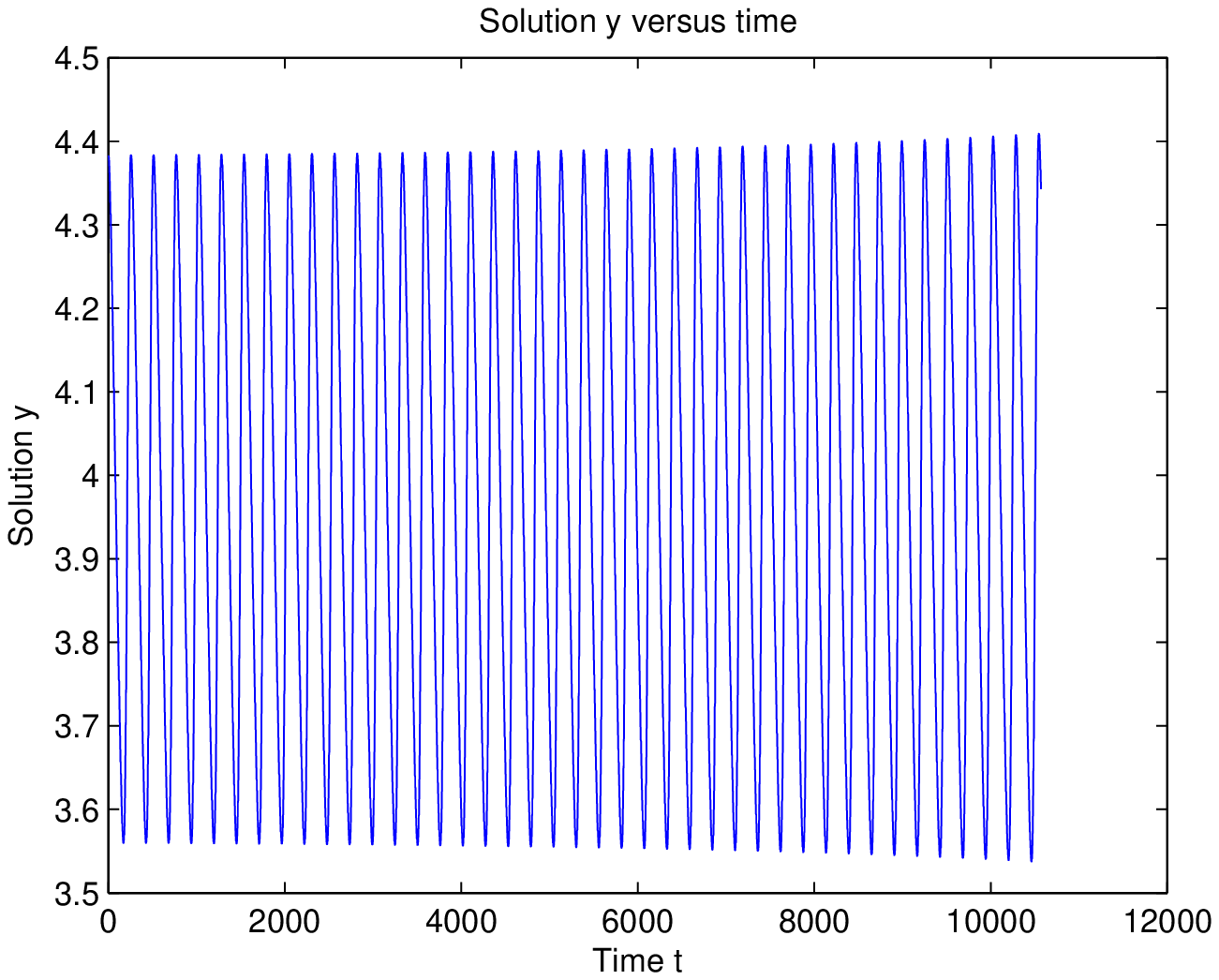}\includegraphics[width=0.49\linewidth]{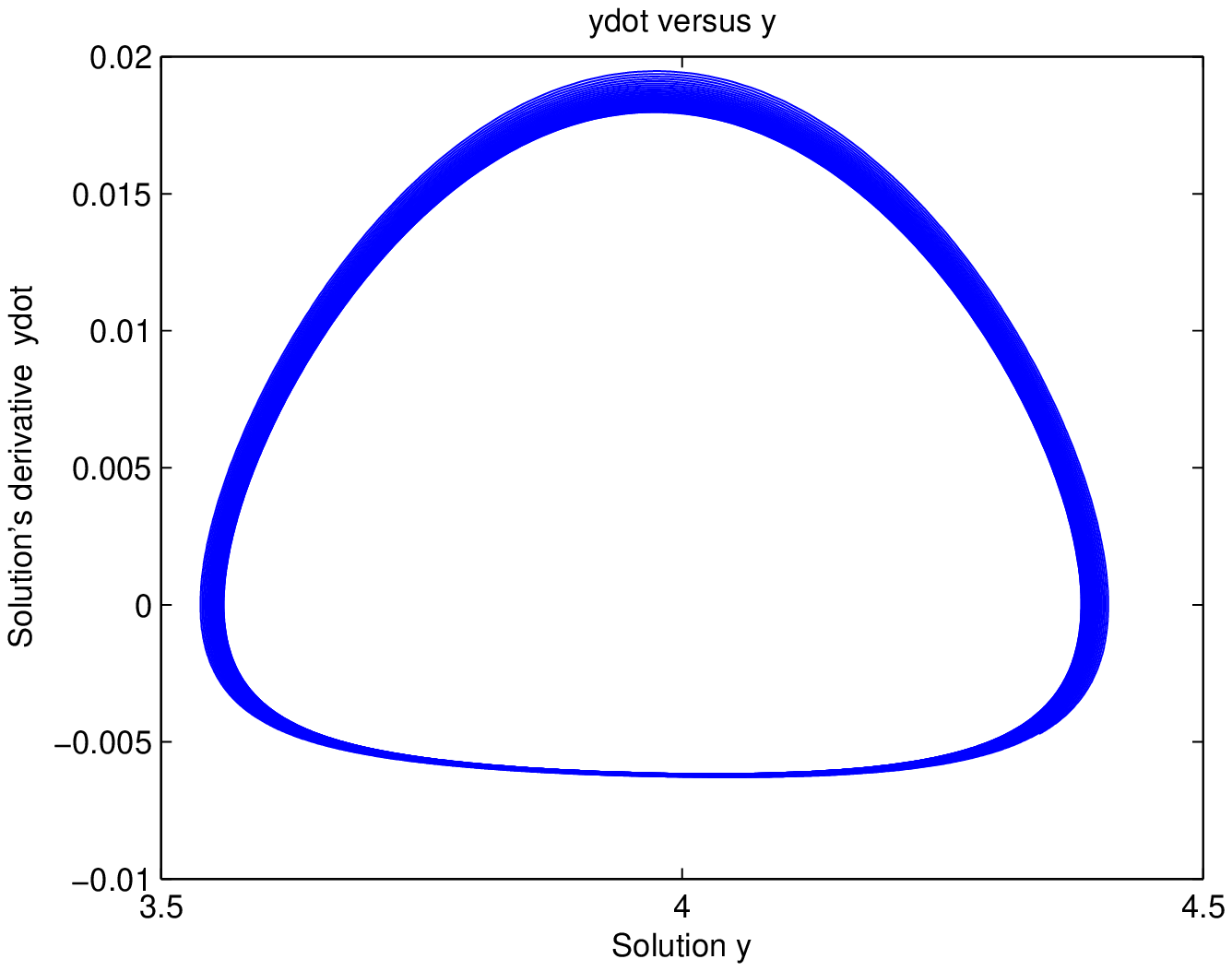}
\caption{Graphs of $x(\cdot),$ respectively $y(\cdot)$, for point
$P_3,\;\;c=0.425$ .}
\end{figure}

\begin{figure}\centering
\includegraphics[width=0.49\linewidth]{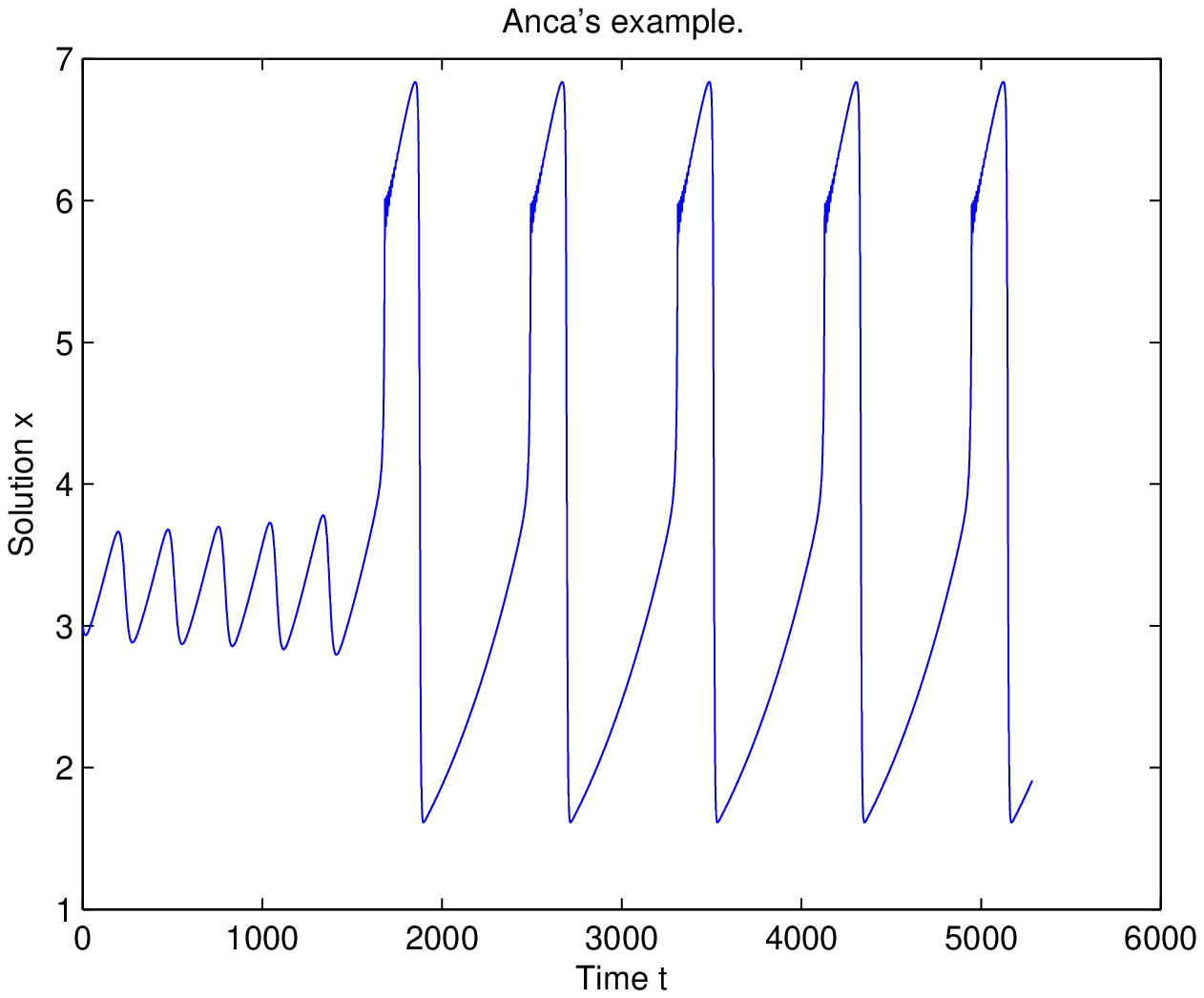}\includegraphics[width=0.49\linewidth]{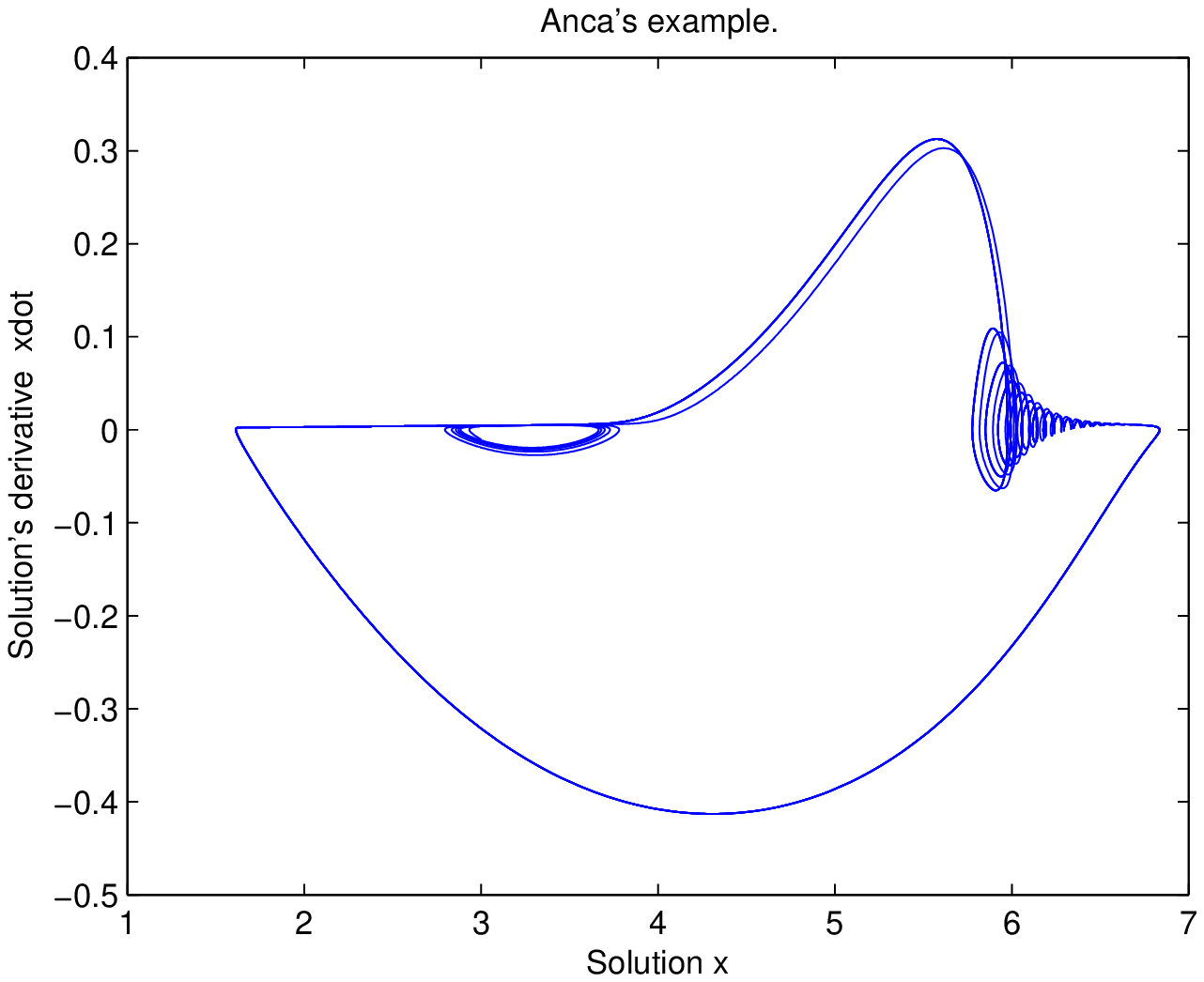}\\
\includegraphics[width=0.49\linewidth]{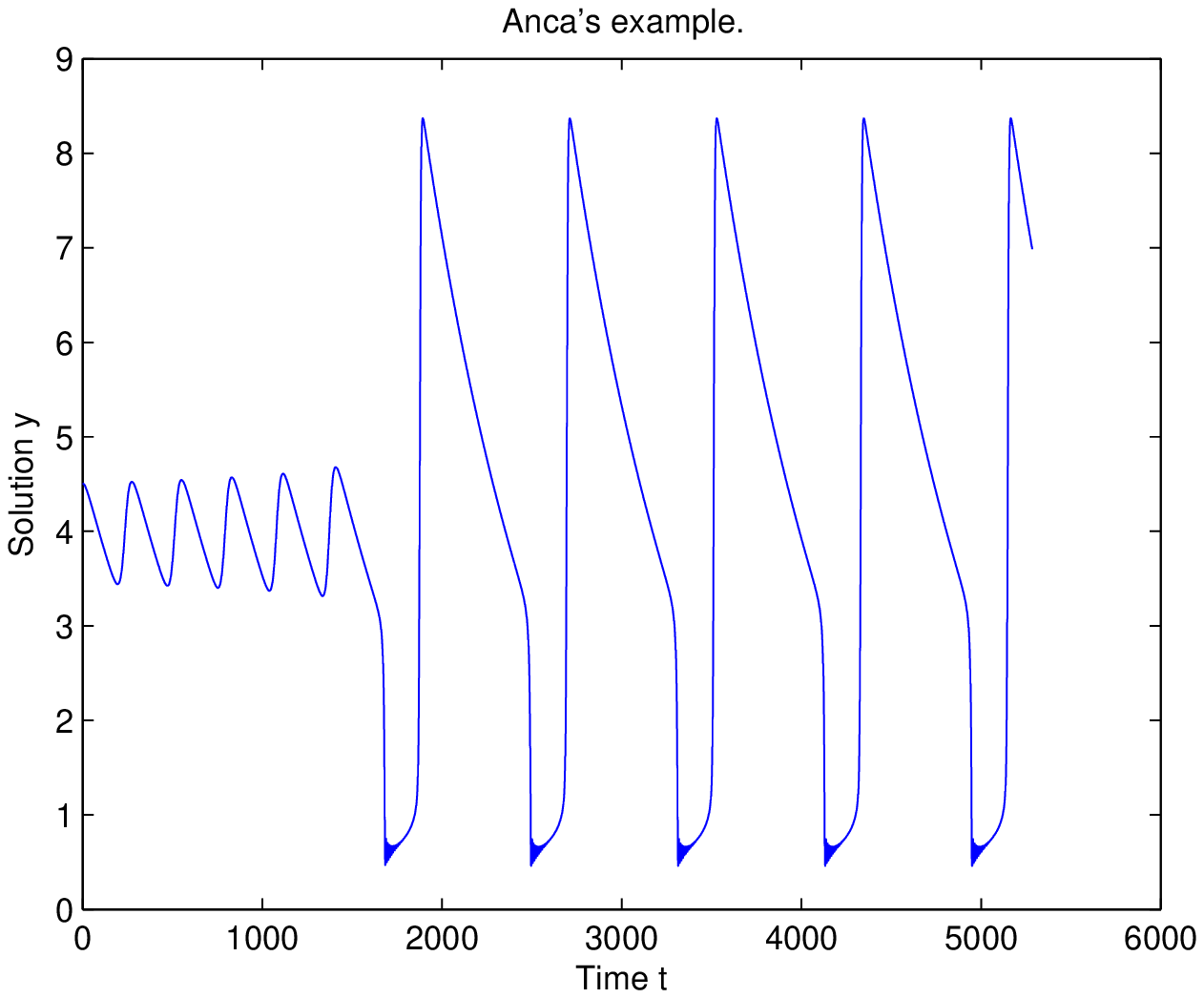}\includegraphics[width=0.49\linewidth]{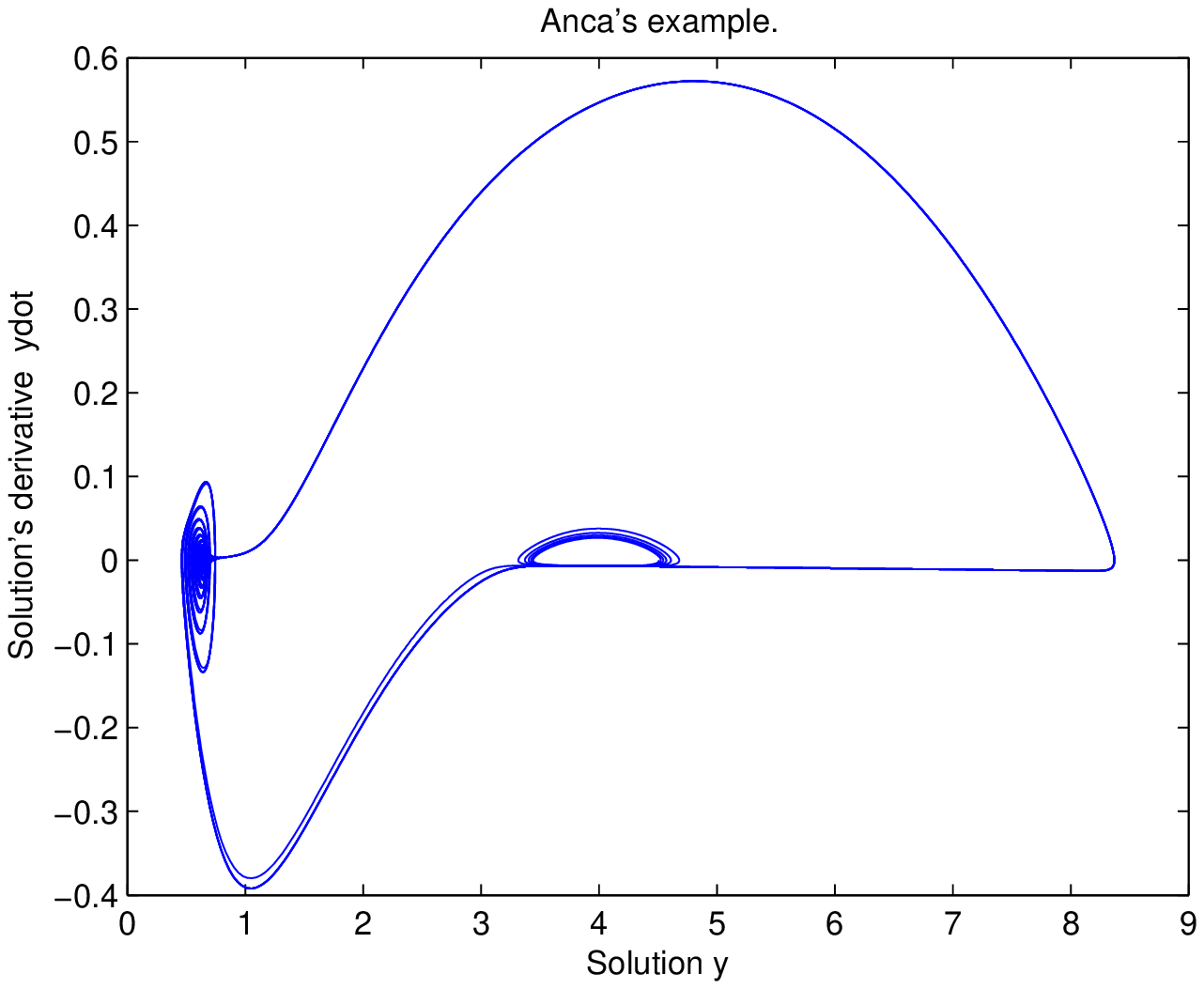}
\caption{Graphs of $x(\cdot),$ respectively $y(\cdot)$, for point
$P_3,\;\;c=0.55$.}
\end{figure}
In the present paper, for the same point, $P^*$,  we show, in
parallel, the behaviour of $y$ and that of $x$ for the parameters
chosen in the (most interesting) zone, corresponding to the zone 3
of the bifurcation diagram (Fig. 6, left). There, for $y$, two limit
cycles exist (one inside the other), the inner cycle being unstable
while the outer one, stable. By numerical investigations, we found
that point $P_3(7.55,\,0.0015)$ lies in this zone. For this point,
the second equilibrium of system \eqref{eq-x}, \eqref{eq-y} is
$x_2=3.24777,\,y_2=3.95811$.\\

As initial function $y_0=\phi$, for eq. \eqref{eq-y}, we considered
functions of the type  $\phi(s)=y_2+ce^{\mu s}\cos(\omega
s),\,s\in[-r,0]$, where $\lambda_{1,2}=\mu\pm i\omega$ are the two
eigenvalues with greatest real part corresponding to the chosen
parameters. Hence the parameter $c$ determines the distance between
the initial function $\phi$ and the equilibrium $y_2.$

Depending on $c$, the orbit $y(t,\phi)$  will either tend towards $y_2$ (being rejected by the unstable inner limit cycle) or towards the outer limit cycle. First we take $c=0.2,$ and $c=0.41$ and see that $y(t)$ tend towards $y_2$ and $x(t)$ tends to $x_2$ (Figs 11 and 12). For $c=0.425$ the oscillating orbits have increasing amplitude (Fig. 13). It is clear that between $c=0.41$ and $c=0.425$ there must be a value $c^*$ such that for the  initial function $\phi^*(s)=y_2+c^*e^{\mu s}\cos(\omega s),\,s\in[-r,0]$, the corresponding solution $y(t,\phi^*)$ of eq. \eqref{eq-y} is a periodic function. This is the unstable periodic orbit (inner limit cycle on the center manifold) of the equation for $y$.

We notice that for $c=0.425$ the amplitude of the orbits increases very slowly. If we take a greater value of $c$, e.g. $c=0.55,$  for a short interval of time, the amplitude increases slowly, and then there is a sudden increase of the amplitude, and the solution rapidly tends to a periodic orbit, that is stable (Fig. 14). This one corresponds to the outer limit cycle from the center manifold.
This behaviour of $y$ is reproduced, qualitatively, by $x$, as is seen in Figure 14.

\section{Conclusions}

In this paper we investigated the behaviour of the solution $x$ of the system \eqref{eq-x}, \eqref{eq-y}, the function $x$ representing the nondimensional density of proliferating cells in the leukemia disease. It is seen that the behaviour of $x$ is determined by that of $y$ (density of resting cells).

We  studied the stability of the equilibrium solutions, and we
presented the conditions that the parameters should satisfy, in
order that the equilibrium solutions $x_1,\,x_2$ of eq. \eqref{eq-x}
be stable or even globally asymptotically stable.

Then, we establish the conditions on the parameters such the solution $x$ mimics a Hopf bifurcation (i.e. the equilibrium loses stability and a stable periodic orbit occurs).

The third type of interesting behaviour of $y$ previously studied \cite{AVI-RMG} is the Bautin bifurcation. We showed that, when $y$ presents some Bautin bifurcation, the solution $x$ also mimics this kind of bifurcation.

The study concerning the Bautin bifurcation is important because, for parameters in the domains corresponding to the domain 3 of the bifurcation diagram, the behaviour of $x$ may have two very different aspects. That is $x(t)$ may tend to $x_2$ or to the periodic orbit that corresponds to the outer limit cycle of the reduced to the center manifold problem for $y$. The initial value of $x$ is irrelevant in deciding what behaviour $x$ has. Only the initial function for $y$ matters.

\vspace{0.3cm}

We also provide tables and figures with points of Bautin bifurcation space for $n=2$ and  $\beta_0\in\{0.5,\,1,\,1.5,\,2\},$  $k\in\{1.1,\,1.2,...,1.9\}$.
The conclusions presented for the occurrence or non-occurrence of Bautin bifurcation, refer to an important zone of the parameters space, that is $1\leq n \leq 12$ and $\beta_0,\,k$ as above.

The investigation may be extended to other parameters zones by using the method exposed in \cite{AVI-RMG}.

More than that, we also performed numerical integrations that confirm the qualitative behaviours that we theoretically predicted - Figs. 2-5 and Figs. 11-14.

\vspace{0.3cm}

\textbf{Remark.} \textit{We remark once more that, for all considered problems: stability, Hopf bifurcation or Bautin bifurcation, the long time behaviour of $x$ does not depend on the initial condition for this function ($x(0)=x_0$), but depends on the initial state of the function $y$.}

\textit{This is important for the physician or for the biologist, because it shows that, in order to make a long time prediction on the development of the illness, the initial value of $y$, that is of the function that models the "resting cells" should be known and not that of $x$.  As an example, in Figures 11-14, for the same parameters, depending on the initial condition of $y$, $x(t)$ tends either to $x_2$ or to a large stable periodic orbit. The initial value of $x(\cdot)$ may only influence the length of the time interval in which $x(t)$ approaches its asymptotic behaviour, determined by $y$.
}

In the Appendix we prove some (simple) results that show how $y$ detemines $x$.

\vspace{0.3cm}

\section{Appendix}

We assert and prove some  results that are necessary in the  study
of the main part of the paper.

In order to present the results in an unified form, we make the
following notations. We denote by $(x^*,\,y^*)$ one of the two
equilibria $(x_1,\,y_1)$ or $(x_2,\,y_2)$. Remind that we denoted
the terms containing $y$ from eq. \eqref{eq-x}, by
$\widetilde{F}(y)(t)$, such that eq. \eqref{eq-x} takes the form
\[\frac{d x}{dt}(t)=-\gamma x(t) +\widetilde{F}(y)(t).
\]

If $y^*= y_2,$ we rewrite eq. \eqref{eq-x} as
\[\frac{d(x-x_2)}{dt}=-\gamma(x-x_2)+\frac{\beta_0y(t)}{1+y(t)^n}-\frac{\beta_0y_2}{1+y_2^n}-
\frac{k}{2}\left(\frac{\beta_0y(t-r)}{1+y(t-r)^n}-\frac{\beta_0y_2}{1+y_2^n}\right).
\]
We set $u=x-x_2,\;\;$  and the equation becomes
\[\frac{d u}{dt}(t)=-\gamma u(t) +\widetilde{F}(y)(t)-\widetilde{F}(y_2).
\]
We denote
\[H(y)(t):=\left\{
          \begin{array}{ll}
            \widetilde{F}(y)(t), & \hbox{if}\,\,\, y^*=y_1, \\
            \widetilde{F}(y)(t)-\widetilde{F}(y_2), & \hbox{if}\,\,\,  y^*=y_2,
          \end{array}
        \right.
\]
and we see that in both situations above our equation has the form
\begin{equation}\label{new-eq-y}\frac{d u}{dt}(t)=-\gamma u(t) +H(y)(t),
\end{equation}
where $H(y)(t)\rightarrow 0$ when $t\rightarrow \infty.$

\vspace{0.3cm}

\textbf{Proposition 2.1.} \textit{If $y(t)\rightarrow y^*$ when $t\rightarrow \infty$, then $x(t)\rightarrow x^*$ when} $t\rightarrow \infty$.

\vspace{0.1cm}

\textbf{Proof.} We prove that $u(t)\rightarrow 0$ when $t
\rightarrow \infty$. The solution, denoted $u(\cdot,\,u_0)$, of the
above equation with the initial condition $u(0)=u_0$, is
\begin{equation}\label{sol-int}u(t,\,u_0)=u_0e^{-\gamma t}+e^{-\gamma t}\int_0^t e^{\gamma s}H(y)(s)ds.
\end{equation}

Let $\varepsilon >0$.  Since $H(y)(t)\rightarrow 0$  there is a $T>0$ such that for $t>T,$ we have $|H(y)(t)|<\frac{\varepsilon\gamma}{3}$. Then we write, for $t \geq T,$
\begin{equation}\label{tends0}|u(t,u_0)|\leq |u_0e^{-\gamma t}|+|e^{-\gamma t}\int_0^T e^{\gamma s}H(y)(s)ds|+|e^{-\gamma t}\int_T^t e^{\gamma s}H(y)(s)ds|\leq
\end{equation}
\[\leq |u_0|e^{-\gamma t}+ e^{-\gamma t}\int_0^T e^{\gamma s}|H(y)(s)|ds+e^{-\gamma t}\int_T^t e^{\gamma s}\frac{\varepsilon\gamma}{3} ds=
\]
\[=|u_0|e^{-\gamma t}+ e^{-\gamma t}\int_0^T e^{\gamma s}|H(y)(s)|ds+\frac{\varepsilon\gamma}{3}
\frac{1}{\gamma}\left(1-e^{\gamma(T-t)} \right).
\]
We chose a $T'>T$ such that $|u_0|e^{-\gamma T'}<\frac{\varepsilon}{3}$ and $e^{-\gamma T'}\int_0^T e^{\gamma s}|H(y)(s)|ds<\frac{\varepsilon}{3}$. \\
Then, for any $t\geq T'$, we have
\[|u(t,u_0)|<\varepsilon,
\]
and the conclusion follows. $\Box$

\vspace{0.5cm}
\textbf{Proposition 2.2.} \emph{If $y(\cdot)$ is periodic, there is a initial value $\widetilde{x}_0$ of eq. \eqref{eq-x} such that the corresponding solution $x(\cdot,\,\widetilde{x}_0)$ is periodic and it is asymptotically stable.}

\vspace{0.2cm}

\textbf{Proof.} If $y(\cdot)$ is periodic, then the function $H(\cdot)$ of eq. \eqref{new-eq-y} is periodic. Then, a classical result (\cite{H}, pg. 214)  shows that there is an $\widetilde{u}_0$ such that the corresponding solution $u(\cdot,\,\widetilde{u}_0)$ of eq. \eqref{new-eq-y} is periodic.

We show that it is asymptotically stable. Let us consider another initial condition $u_0.$ The corresponding solution of \eqref{new-eq-y} is
\[u(t,u_0)=e^{-\gamma t}u_0+\int_0^te^{\gamma (s-t)}H(y)(s)ds=
 \]
\[=e^{-\gamma t}\widetilde{u}_0+\int_0^te^{\gamma (s-t)}H(y)(s)ds+e^{-\gamma t}(u_0-\widetilde{u}_0)=
\]
\[=u(t,\,\widetilde{u}_0)+e^{-\gamma t}(u_0-\widetilde{u}_0),
\]
since the sum of the first two terms is the solution of \eqref{new-eq-y} corresponding to $\widetilde{u}_0$. It is clear that $$\lim_{t\rightarrow \infty}|u(t,u_0)-u(t,\,\widetilde{u}_0)|=0, $$
that implies the asymptotic stability of the periodic solution.

The conclusions upon $u$ transfer to $x,$ q.e.d. $\Box$

\vspace{0.5cm}

\textbf{Proposition 2.3.} \textit{If $y(\cdot)$ tends to a periodic function when $t\rightarrow \infty$, then $x(t,y_0)$ tends to a periodic orbit when $t\rightarrow \infty$.
}

\vspace{0.3cm}

\textbf{Proof.} Let $p:[0,\infty)\mapsto \mathbb{R} $ be a periodic function  such that $$\displaystyle\lim_{t\rightarrow \infty}|y(t)-p(t)|=0.$$

With the notation of the preceding Proposition, we see that $H(p)(\cdot)$ is periodic.
We set  $o(t)=H(y)(t)-H(p)(t)$.

It is not difficult to show that $\displaystyle\lim_{t\rightarrow \infty}o(t)=0.$

Now we write eq. \eqref{new-eq-y} as
 $$\dot{u}(t)=-\gamma u(t) + H(p)(t)+o(t) $$ that has solution
 $$u(t,u_0)=e^{-\gamma t}u_0+\int_0^te^{\gamma (s-t)}H(p)(s)ds+\int_0^te^{\gamma (s-t)}o(s)ds. $$

As we pointed out in the proof of Proposition 2.2, in \cite{H} it is shown that there is an initial value $\widetilde{u}_0$ such that the function
\[e^{-\gamma t}\widetilde{u}_0+\int_0^te^{\gamma (s-t)}H(p)(s)ds
\]
is periodic. On the other hand, by a reasoning similar to that applied to the last term of \eqref{tends0}, we prove that $\displaystyle\lim_{t\rightarrow \infty}\int_0^te^{\gamma (s-t)}o(s)ds=0.$

It follows that for the initial condition $\widetilde{u}_0,$ the solution $u(\cdot, \widetilde{u}_0)$ tends to a periodic function.
Now, for $u_0\neq \widetilde{u}_0,$ we have
\[u(t)=e^{-\gamma t}u_0+\int_0^te^{\gamma (s-t)}H(p)(s)ds+\int_0^te^{\gamma (s-t)}o(s)ds=
 \]
\[=e^{-\gamma t}\widetilde{u}_0+\int_0^te^{\gamma (s-t)}H(p)(s)ds+e^{-\gamma t}(u_0-\widetilde{u}_0)+\int_0^te^{\gamma (s-t)}o(s)ds.
\]
Since the last two terms tend to zero when $t \rightarrow \infty,$ it again follows that $u(\cdot,u_0)$ tends to a periodic function. Hence $x(\cdot,x_0)$ has the same property.$\Box$

\end{document}